\documentclass[final]{amsart}
\usepackage[T1]{fontenc}
\usepackage{amsmath}
\usepackage{amsthm}
\usepackage{amssymb}
\usepackage{mathtools}
\usepackage{accents}
\usepackage{cases}
\usepackage{empheq}
\usepackage{pgfplots}
\usepgfplotslibrary{groupplots}
\usepgfplotslibrary{fillbetween}
\pgfplotsset{compat=1.14}
\usepackage{fullpage}
\usepackage{subcaption}
\usepackage{booktabs}
\usepackage{fullpage}
\usepackage{graphicx}
\usepackage{multirow}

\usepackage[dvipsnames]{xcolor}
\usepackage[colorlinks=true, linkcolor=Purple, citecolor=Orange, urlcolor=Purple]{hyperref}
\mathtoolsset{showonlyrefs}

\newcounter{custom}
\numberwithin{custom}{section} 

\let\oldsubsection\subsection
\renewcommand{\subsection}{\stepcounter{custom}\oldsubsection}

\newtheorem{remark}[custom]{Remark}

\usepackage{tkz-graph}
\usepackage{tikz}
\usepackage{xcolor}
\usetikzlibrary{arrows}
\usepackage[ruled,linesnumbered]{algorithm2e}
\usepackage[noend]{algpseudocode}
\usepackage{wrapfig}
\usepackage{caption}

\usepackage{tristan}
\usepackage{tikz-3dplot}
\usetikzlibrary{matrix}
\usepackage{todonotes}

\definecolor{plot1}{HTML}{05A8AA}
\definecolor{plot2}{HTML}{B8D5B8}
\definecolor{plot3}{HTML}{D7B49E}
\definecolor{plot4}{HTML}{DC602E}
\definecolor{plot5}{HTML}{BC412B}

\SetKwComment{Comment}{$\triangleright$\ }{}

\renewcommand{\vec}[1]{\ensuremath{\boldsymbol{#1}}}

\title[]{Diffusion Synthetic Acceleration for polytopal discretisations of Boltzmann transport}
\author[Calloo, A. et al.]{Ansar Calloo$^{1}$, Matthew Evans $^{2*}$, Fran\c{c}ois Madiot$^{3}$, Tristan Pryer$^{2,4}$}
\address{$^1$Universit\'e Paris-Saclay, CEA, Service de G\'enie Logiciel pour la Simulation, 91191, Gif-sur-Yvette, France}
\address{$^2$Mathematical Sciences, University of Bath}
\address{$^3$Universit\'e Paris-Saclay, CEA, Service d'\'Etudes des R\'eacteurs et de Math\'ematiques Appliqu\'ees, 91191, Gif-sur-Yvette, France}
\address{$^4$Institute of Mathematical Innovation, University of Bath}
\email{$^*$\url{mje45@bath.ac.uk}}
\date{}

\begin{document}

\bibliographystyle{acm}

\begin{abstract}
  We present a computational study of diffusion synthetic acceleration
  (DSA) for the monoenergetic, isotropically scattering $S_N$
  transport equations, discretised in space by a polytopic
  discontinuous Galerkin method. Using a discrete ordinates angular
  discretisation, we construct the DSA correction with an
  interior-penalty diffusion operator and compare a classical
  symmetric interior penalty (SIP) formulation with a modified
  interior penalty (MIP) variant, together with homogeneous Dirichlet
  and Marshak (Robin) diffusion boundary conditions imposed weakly in
  the DG framework.

  We quantify the observed convergence behaviour of the resulting
  source iteration across variations in optical thickness, scattering
  ratio, angular quadrature, mesh refinement, polynomial degree and
  mesh anisotropy on families of bounded Voronoi meshes. The results
  show that MIP-based DSA remains robust across the parameter ranges
  tested, whereas SIP-based DSA can lose robustness in the
  intermediate regime. In challenging optically thick, highly
  scattering settings, the observed convergence factors for the
  MIP-based schemes are typically below $0.6$.
\end{abstract}

\maketitle

\section{Introduction}

Transport models based on the linear Boltzmann transport equation
(BTE) sit at the core of three application areas that drive both the
modelling choices and the computational constraints in this paper,
nuclear engineering (reactor physics, criticality support calculations
and shielding), medical physics (in particular radiotherapy dose and
optimisation workflows) and space technologies (space radiation
transport for shielding design and crew risk
assessment). Deterministic discrete ordinates ($S_N$) solvers remain a
mainstay in these settings, both as production tools and as components
inside hybrid workflows, for example in large-scale shielding
calculations and reactor analysis \cite{evans2010denovo}. In
radiotherapy, there is renewed interest in fast deterministic
transport surrogates for dose calculation and optimisation loops,
including proton-specific transport modelling pipelines
\cite{ashby2025efficient}.

The discrete ordinates method (DOM) is a standard approach to
semi-discretise the BTE by approximating its angular dependence with a
quadrature rule, yielding a coupled system of spatial equations (the
$S_N$ equations). Given the scale of these systems, iterative methods
with preconditioning are required. One widely adopted scheme is source
iteration, a Richardson iteration that avoids full inversion by
isolating the transport operator, which can be inverted efficiently by
a sweep, either directly or in a matrix-free manner. Source iteration
performs well in many settings but deteriorates in diffusive regimes
such as moderators, reflectors and shielding problems. In highly
scattering and optically thick settings the contraction factor
approaches unity, leading to slow convergence or
pseudo-convergence. Under standard absorption assumptions, source
iteration is contractive, but in highly scattering and optically thick
regimes its convergence factor can be arbitrarily close to one, making
acceleration essential \cite{dsa_lifeline}.

Acceleration strategies for the BTE are often grouped into synthetic
methods, which preserve the discrete fixed point through an additive
correction, and other acceleration strategies, which may improve
convergence without corresponding exactly to a fixed-point-preserving
preconditioner. Examples include transport synthetic acceleration
\cite{ramone1997transport} and quasi-diffusion acceleration
\cite{gol1964quasi}. Synthetic methods approximate the inverse of the
full transport system and yield a preconditioned fixed-point iteration
via an additive correction. A particularly effective option in
diffusive regimes is diffusion synthetic acceleration (DSA)
\cite{alcouffe1977diffusion}. DSA can deliver mesh- and
thickness-robust convergence, but only when the diffusion correction
is discretised in a way that is compatible with the transport
discretisation.

For discontinuous Galerkin (DG) transport discretisations with upwind
numerical fluxes, the dominant discrete dissipation in thick regimes
is delivered facewise by the upwind jump terms. A diffusion correction
based on a symmetric interior penalty (SIP) DG discretisation is
elliptically coercive, but elliptic coercivity alone does not
guarantee effective acceleration. On some meshes or for some
quadrature--mesh alignments the SIP jump control can be weaker than
the dissipation delivered by the upwind transport flux, which degrades
or stalls or even destabilises the accelerated iteration. This
motivates penalty choices that explicitly match, face by face, the
transport dissipation scale while retaining elliptic stability of the
diffusion solve. A companion analytic paper \cite{dsaanal} develops
this perspective and provides contraction estimates for
transport-matched penalty choices.

Polytopic DGFEMs are attractive for transport because they support
high-order approximation, flexible treatment of complex geometries and
a rigorous approximation theory \cite{cangiani2014hp}. Although DG
methods introduce additional face terms, the ability to work on
agglomerated or general polytopic meshes (for example, Voronoi
tessellations) can reduce element counts in smooth regions or
geometrically complex settings, while retaining parallel scalability
in transport sweeps. Recent work has also clarified when such meshes
admit robust sweep constructions \cite{calloo2025cycle}. Despite these
advantages, DSA-based preconditioning for polytopic DGFEMs remains
relatively unexplored, particularly from the standpoint of quantifying
how penalty choices, $hp$-refinement, quadrature and mesh anisotropy
affect the contraction factor.

This computational study explores diffusion synthetic acceleration for
the monoenergetic, isotropically scattering $S_N$ transport equations
discretised in space by a polytopic discontinuous Galerkin method. The
present paper accompanies the companion paper \cite{dsaanal}, which
focusses on the analytic perspective, establishes discrete contraction
bounds and motivates transport-matched penalty choices.

We implement DSA in the same discontinuous polynomial space and on the
same mesh as the transport discretisation, and we focus on
interior-penalty diffusion operators with weakly imposed boundary
conditions. In particular, we compare a classical symmetric interior
penalty (SIP) diffusion discretisation with a modified interior
penalty (MIP) variant. The MIP idea is well established in the DSA
literature for DG $S_N$ discretisations on simplicial and Cartesian
meshes \cite{Wang01102010,ragusa2010two}, where it is used to prevent
loss of robustness in optically thick and highly scattering
regimes. Our focus is not novelty of MIP per se, but rather a
systematic investigation of its behaviour in a modern polytopic DG
setting, where facet geometry, local anisotropy and quadrature--mesh
interactions can substantially affect penalty scaling and sweep
performance. To the best of our knowledge, this has not been explored
in a dedicated empirical study on general polytopic meshes.

We also examine the role of the diffusion boundary model by comparing
a strong vacuum surrogate (homogeneous Dirichlet) with a Marshak
(Robin) condition motivated by the diffusion limit of vacuum inflow.

We quantify how the empirical contraction factor depends on the key
parameters that control difficulty in the diffusive regime. This
includes the angular quadrature (choice of $S_N$ rule), polynomial
degree $p$, scattering ratio and mesh refinement on families of
bounded Voronoi meshes. We also study the effect of element
anisotropy, including regimes in which interior-penalty schemes can
under-penalise because facets are small relative to element
diameter. We report its impact on observed convergence factors,
iteration counts and wall-clock time. We report the impact of aspect
ratio on the observed contraction factor and, within a fixed
direct-solve implementation, on practical performance metrics
including iteration counts and wall-clock time. To isolate
outer-iteration behaviour from differences in inner linear-solver
technology, we use direct sparse solves for both transport and
diffusion throughout the reported experiments.

DG spatial discretisations have been used for discrete ordinates
($S_N$) transport since the earliest triangular-mesh formulations,
with supporting analysis developed soon afterwards
\cite{reed1973triangular,lesaint1974finite}. Subsequent work extended
these ideas to three-dimensional unstructured settings and to element
families commonly used in reactor and shielding calculations,
including tetrahedral and hexahedral decompositions, with further
study of stability and convergence behaviour in practical
configurations \cite{wareing2001discontinuous,Wang01092009}.

Acceleration has an equally long history in neutronics. DSA was
introduced as a synthetic method for accelerating source iteration in
scattering-dominated regimes \cite{alcouffe1977diffusion} and its
behaviour in the diffusive limit, including the role of discretisation
consistency, has been studied extensively
\cite{larsen1982unconditionally,adams1992diffusion,morel2005sn,dsa_lifeline}. A
key theme in this body of work is that robust acceleration in
optically thick, highly scattering regimes requires the diffusion
correction to control the same slow-to-damp error components as the
underlying transport iteration, while preserving the correct discrete
fixed point.

For DG $S_N$ discretisations, several diffusion-correction constructions
have been proposed that aim to remain stable and effective across both
transport-dominated and diffusion-dominated regimes. Discontinuous
diffusion operators based on interior penalty ideas have been developed
for high-order DG $S_N$ transport on locally refined unstructured meshes,
with particular emphasis on penalty choices that avoid loss of robustness
in thick regimes \cite{Wang01102010,ragusa2010two}. In parallel, iterative
solution strategies for the large coupled systems arising from high-order
and $hp$ DG Boltzmann discretisations have been investigated, including
fixed-point and Krylov-based approaches and their interaction with
physics-informed preconditioning \cite{houston2024iterative}. The impact
of material heterogeneity and discontinuities on performance and
robustness has also been studied, including variants designed to remain
effective in heterogeneous configurations
\cite{warsa2004krylov,Southworth01022021}. These developments connect
naturally with scalable solvers for the auxiliary diffusion problem,
including preconditioned conjugate gradient and multigrid strategies,
which are often essential in large-scale applications
\cite{hackemack2016higher}.

More recent work has focused on extending DG transport discretisations
to general polytopic meshes, motivated by geometric flexibility and by
practical mesh constructions such as Voronoi tessellations
\cite{TURCKSIN2014356,hackemack2018quadratic,cangiani2014hp}. In this
direction, quadrature-free and general polytopic DG formulations have
been developed for transport problems, helping to reduce geometry-driven
assembly costs and to broaden the class of admissible cell shapes
\cite{radley2024quadrature}. Unified high-order space--angle--energy
polytopic DG discretisations for linear Boltzmann transport have also
been developed, providing a modern setting in which sweep-based transport
solves and acceleration technology can be combined on general polytopic
meshes \cite{houston2024efficient}. While DSA constructions tailored to
polygonal meshes have been proposed \cite{TURCKSIN2014356}, systematic
empirical studies that isolate how contraction depends on penalty choice,
boundary modelling, $hp$-refinement, quadrature and mesh anisotropy in a
fully polytopic DG transport setting remain limited.

The present work sits at this intersection. We build computational
tests that probe precisely the mechanisms that are known to govern
stability and effectiveness of interior-penalty DSA, including
robustness of penalty scalings in adverse diffusion settings
\cite{dong2022robust}, and we interpret the observed contraction
behaviour through the discrete contraction perspective developed in
\cite{dsaanal}.

The rest of the paper is organised as follows: \S
\ref{section:transport_dg} presents the polytopic DGFEM discretisation
of the $S_N$ transport equations and the associated source iteration
scheme. \S \ref{section:dsa_dg} introduces the DSA diffusion
correction, including the interior-penalty discretisation, diffusion
boundary conditions and the SIP/MIP penalty choices used throughout,
with implementation remarks guided by \cite{dsaanal}. Numerical
experiments are reported in \S \ref{section:numerics}, where we
quantify empirical contraction factors. \S\ref{section:conclusion}
concludes.

\section{Discontinuous Galerkin scheme for the transport equation}
\label{section:transport_dg}

We begin by recalling the steady-state, monoenergetic linear Boltzmann
transport equation with isotropic scattering. Let
$\Omega\subset\mathbb{R}^d$, $d\in\{2,3\}$, be a bounded Lipschitz
domain with boundary $\Gamma:=\partial\Omega$ and outward unit normal
$\vec{n}(\vec{x})$. Let $\sigma_t$ and $\sigma_s$ denote the
macroscopic total and scattering cross-sections. We assume that the
material coefficients are independent of direction, so
$\sigma_t,\sigma_s$ depend only on $\vec{x}\in\Omega$, and that
scattering is isotropic, meaning the scattering kernel is constant.
Accordingly, we take $\sigma_t,\sigma_s\in L^\infty(\Omega)$ and
assume the standard bounds
\begin{equation}
    \label{eqn:xs_bounds}
    \sigma_t(\vec{x})\ge \sigma_{t,\min}>0,\qquad 0\le \sigma_s(\vec{x})\le \sigma_t(\vec{x}) \quad \text{for a.e. }\vec{x}\in\Omega.
\end{equation}
We define the absorption cross-section by $\sigma_a :=
\sigma_t-\sigma_s\ge 0$.

For an angular flux $\psi(\vec{x},\vec{\omega})$ we define the scalar
flux as the uniform angular average
\[
\phi(\vec{x}) := \langle \psi(\vec{x},\cdot)\rangle
:= \frac{1}{|\mathbb{S}^{d-1}|}\int_{\mathbb{S}^{d-1}}\psi(\vec{x},\vec{\omega}') d\vec{\omega}'.
\]
Given a source $f\in L^2(\Omega\times\mathbb{S}^{d-1})$ and inflow
boundary data $g_{\mathrm{D}}$ (in the numerical experiments we
typically take vacuum inflow $g_{\mathrm{D}}\equiv 0$), we seek $\psi$
satisfying
\begin{equation}
    \label{eqn:bte}
    \begin{split}
        \vec{\omega}\cdot\nabla\psi(\vec{x},\vec{\omega}) + \sigma_t(\vec{x})\psi(\vec{x},\vec{\omega})
        &= \sigma_s(\vec{x})\phi(\vec{x}) + f(\vec{x},\vec{\omega})
        \qquad \text{in }\Omega\times\mathbb{S}^{d-1},\\
        \psi(\vec{x},\vec{\omega}) &= g_{\mathrm{D}}(\vec{x},\vec{\omega})
        \qquad \text{on }\Gamma_-,
    \end{split}
\end{equation}
where the inflow boundary is the phase-space set
\[
\Gamma_- := \{(\vec{x},\vec{\omega})\in\partial\Omega\times\mathbb{S}^{d-1} : \vec{\omega}\cdot\vec{n}(\vec{x}) < 0\}.
\]
Under \eqref{eqn:xs_bounds} and sufficient data regularity, the
transport problem \eqref{eqn:bte} is well posed
\cite{calloo2025cycle}.

\subsection{$S_N$ transport equations}

We discretise \eqref{eqn:bte} in angle using the discrete ordinates
method (DOM). For a given $N_Q\in\mathbb{N}$, let
$\{(w_k,\vec{\omega}_k)\}_{k=1}^{N_Q}$ be a quadrature rule on
$\mathbb{S}^{d-1}$ with normalised weights, meaning
\begin{equation}
    \label{eqn:quadrature_rule}
    \frac{1}{|\mathbb{S}^{d-1}|}\int_{\mathbb{S}^{d-1}} z(\vec{\omega}) d\vec{\omega}
    \approx \sum_{k=1}^{N_Q} w_k  z(\vec{\omega}_k),
    \qquad \sum_{k=1}^{N_Q} w_k = 1,
\end{equation}
for sufficiently regular $z$. Equivalently, $\int_{\mathbb{S}^{d-1}}
z(\vec{\omega}) d\vec{\omega}\approx
|\mathbb{S}^{d-1}|\sum_{k=1}^{N_Q} w_k z(\vec{\omega}_k)$. Commonly
used examples include Chebyshev or trapezoidal quadratures in two
dimensions and level-symmetric quadrature sets in three dimensions. We
refer to \cite{RUKOLAINE2001257,thurgood1995tn} for reviews of
symmetry and exactness properties required to preserve the invariances
of the $(d-1)$-sphere.

For each ordinate direction $\vec{\omega}_k$, we partition the
physical boundary $\Gamma$ into inflow and outflow portions
\begin{equation}
    \Gamma_-^{(k)} := \{\vec{x}\in\Gamma : \vec{\omega}_k\cdot\vec{n}(\vec{x})<0\},
    \qquad
    \Gamma_+^{(k)} := \{\vec{x}\in\Gamma : \vec{\omega}_k\cdot\vec{n}(\vec{x})\ge 0\},
\end{equation}
respectively. Only $\Gamma_-^{(k)}$ enters the DG formulation through
the upwind numerical flux.

Denoting $\psi_k(\vec{x})\approx \psi(\vec{x},\vec{\omega}_k)$ and
$f_k(\vec{x}) := f(\vec{x},\vec{\omega}_k)$, the DOM approximation of
\eqref{eqn:bte} is the coupled system: for $k=1,\dots,N_Q$, find
$\psi_k$ such that
\begin{equation}
    \label{eqn:dom_equations}
    \begin{split}
        \vec{\omega}_k\cdot\nabla\psi_k(\vec{x}) + \sigma_t(\vec{x})\psi_k(\vec{x})
        &= \sigma_s(\vec{x}) \sum_{l=1}^{N_Q} w_l \psi_l(\vec{x}) + f_k(\vec{x})
        \qquad \vec{x}\in\Omega,\\
        \psi_k(\vec{x}) &= g_{\mathrm{D}}(\vec{x},\vec{\omega}_k)
        \qquad \vec{x}\in\Gamma_-^{(k)}.
    \end{split}
\end{equation}
We also define the (semi-discrete) scalar flux by the quadrature average
\begin{equation}
    \label{eqn:scalar_flux}
    \phi(\vec{x}) := \sum_{k=1}^{N_Q} w_k \psi_k(\vec{x})
    \approx \frac{1}{|\mathbb{S}^{d-1}|}\int_{\mathbb{S}^{d-1}} \psi(\vec{x},\vec{\omega}) d\vec{\omega},
\end{equation}
so that the scattering source in \eqref{eqn:dom_equations} can be
written compactly as $\sigma_s(\vec{x}) \phi(\vec{x})$.

\subsection{Discontinuous Galerkin discretisation}

To apply a discontinuous Galerkin discretisation to
\eqref{eqn:dom_equations} we introduce the mesh, trace notation and
the upwind numerical flux.

Let $\mathcal{T}_h$ be a subdivision of $\Omega$ into disjoint open
polytopic elements $\kappa$ such that
$\overline{\Omega}=\cup_{\kappa\in\mathcal{T}_h}\overline{\kappa}$. For
$\kappa\in\mathcal{T}_h$ let $h_\kappa:=\mathrm{diam}(\kappa)$ and
define the global meshsize $h:=\max_{\kappa\in\mathcal{T}_h}
h_\kappa$. We assume $\mathcal{T}_h$ is shape-regular in the sense
that there exists $C_r>0$ such that
\begin{equation}
    \frac{h_\kappa}{\eta_\kappa}\le C_r,\qquad \forall \kappa\in\mathcal{T}_h,
\end{equation}
where $\eta_\kappa$ denotes the diameter of the largest ball contained
in $\kappa$. This assumption does not prevent the presence of facets
whose diameter is arbitrarily small compared with $h_\kappa$, a
feature that is relevant for polytopic meshes and agglomerated
constructions. In addition, we assume the standard local covering
assumptions used in polytopal DG analysis, so that the usual trace,
inverse and approximation estimates hold \cite{cangiani2022}.

On $\mathcal{T}_h$ we define the broken Sobolev space
\begin{equation}
    H^1(\mathcal{T}_h):=\{v\in L^2(\Omega): v|_\kappa\in H^1(\kappa)\ \ \forall\kappa\in\mathcal{T}_h\}.
\end{equation}
For each element $\kappa$ we denote its boundary by $\partial\kappa$
and write $\vec{n}_\kappa(\vec{x})$ for the outward unit normal on
$\partial\kappa$. Given an ordinate direction $\vec{\omega}_k$, we
define the element inflow boundary
\begin{equation}
    \partial_-^{(k)}\kappa := \{\vec{x}\in\partial\kappa : \vec{\omega}_k\cdot\vec{n}_\kappa(\vec{x})<0\}.
\end{equation}

For $v\in H^1(\mathcal{T}_h)$ we denote by $v^+_\kappa$ the interior
trace from $\kappa$ onto $\partial\kappa$. For almost every
$\vec{x}\in\partial\kappa\setminus\Gamma$ there exists a unique
neighbouring element $\kappa'\in\mathcal{T}_h$ such that
$\vec{x}\in\partial\kappa'$, and we define the exterior trace relative
to $\kappa$ by $v^-_\kappa := v^+_{\kappa'}$ on
$\partial\kappa\setminus\Gamma$. When the element is clear from the
context we omit the subscript.

We now introduce the upwind jump of a scalar-valued function $v$ on
the inflow portion of $\partial\kappa$ away from the physical
boundary. For $\vec{x}\in\partial_-^{(k)}\kappa\setminus\Gamma$ we set
\begin{equation}
    \label{eqn:upwind_jump}
    \ujump{v}_k := v^+ - v^-.
\end{equation}

We define the discontinuous finite element space by
\begin{equation}
    \mathbb{V}_{h,p} := \{v_h\in L^2(\Omega): v_h|_\kappa\in\mathbb{P}_{p_\kappa}(\kappa),\ \kappa\in\mathcal{T}_h\},
\end{equation}
where $p_\kappa\in\mathbb{N}_0$. 

For the $S_N$ system we introduce the vector of discrete directional fluxes
\begin{equation}
    \label{eqn:vector_flux}
    \Psi_h := \begin{bmatrix}\psi_{1,h},\hdots,\psi_{N_Q,h}\end{bmatrix}^{T}.
\end{equation}
Using the normalised quadrature weights from
\eqref{eqn:quadrature_rule} (so $\sum_{l=1}^{N_Q} w_l=1$), we define
the scattering coupling functional
\begin{equation}
    s_{h}(\Psi_h,v_h) := \int_\Omega \sigma_s(\vec{x})\left(\sum_{l=1}^{N_Q} w_l \psi_{l,h}(\vec{x})\right) v_h(\vec{x}) d\vec{x}.
\end{equation}
For each $k=1,\dots,N_Q$ we then seek $\psi_{k,h}\in\mathbb{V}_p$ such that
\begin{equation}
    \label{eqn:bilinear_scheme}
    a_{k,h}(\psi_{k,h},v_h) = s_{h}(\Psi_h,v_h) + \ell_k(v_h)\qquad \forall v_h\in\mathbb{V}_p,
\end{equation}
where, following standard upwind DG formulations for transport
\cite{calloo2025cycle,Southworth01022021,wang2016discrete,johnson1983convergence},
the bilinear form and linear functional are
\begin{equation}
    \label{eqn:transport_bilinear}
    \begin{split}
        a_{k,h}(u_h,v_h)
        := \sum_{\kappa\in\mathcal{T}_h}\Bigg(
        &\int_{\kappa} \big(\vec{\omega}_k\cdot\nabla u_h + \sigma_t u_h\big) v_h d\vec{x}
        - \int_{\partial_-^{(k)}\kappa\setminus\Gamma} \big(\vec{\omega}_k\cdot\vec{n}_\kappa\big) \ujump{u_h}_k v_h ds \\
        &\hspace{2.8em}
        - \int_{\partial_-^{(k)}\kappa\cap\Gamma_-^{(k)}} \big(\vec{\omega}_k\cdot\vec{n}\big) u_h v_h ds
        \Bigg),\\
        \ell_k(v_h)
        := \int_\Omega f_k v_h d\vec{x}
        &- \int_{\Gamma_-^{(k)}} \big(\vec{\omega}_k\cdot\vec{n}\big) g_{\mathrm{D}}(\cdot,\vec{\omega}_k) v_h ds.
    \end{split}
\end{equation}
In the implementation, $f_k$ is represented in $\mathbb{V}_p$ by
interpolation onto the local polynomial basis (and likewise for
$g_{\mathrm{D}}(\cdot,\vec{\omega}_k)$ on inflow facets); for brevity
we retain the notation $f_k$ and $g_{\mathrm{D}}$ in
\eqref{eqn:transport_bilinear}.

This formulation yields standard a priori approximation estimates for
each ordinate
\cite{bey1996hp,johnson1986analysis,lesaint1974finite,cangiani_dg},
which can be combined to analyse the fully discrete angular system
\cite{wang2016discrete,dsaanal}. If reflective or albedo boundary
conditions are imposed, they can be incorporated through modified
boundary fluxes in the same weak sense; we return to boundary
modelling later when discretising the diffusion acceleration equation.

\subsection{Source iteration}

We now introduce an algebraic representation of
\eqref{eqn:bilinear_scheme} and the standard source iteration
(Richardson) solver. For each ordinate $k$ let $\vec{\Psi}_{k,h}$
denote the coefficient vector of $\psi_{k,h}\in\mathbb{V}_p$ in the
chosen DG basis (ordered by element and then by local basis index). We
define the discrete scalar flux coefficient vector by
\begin{equation}
    \label{eqn:phi_vec}
    \vec{\Phi}_h := \sum_{m=1}^{N_Q} w_m \vec{\Psi}_{m,h}.
\end{equation}
With this notation, the DG system \eqref{eqn:bilinear_scheme} for each
direction can be written as
\begin{equation}
    \label{eqn:si_full}
    \vec{A}_k \vec{\Psi}_{k,h} = \vec{S} \vec{\Phi}_h + \vec{L}_k,\qquad k=1,\dots,N_Q,
\end{equation}
where $\vec{A}_k$ is the (generally non-symmetric) transport matrix
associated with $a_{k,h}(\cdot,\cdot)$, $\vec{L}_k$ is the load vector
from $\ell_k(\cdot)$, and $\vec{S}$ is the (spatial) mass matrix
associated with the reaction term $\sigma_s(\cdot)$, that is,
$(\vec{S}\vec{\Phi}_h,\vec{v})$ represents $\int_\Omega \sigma_s
\phi_h v_h d\vec{x}$ in coefficient form. In the isotropic-scattering
setting considered here, $\vec{S}$ couples only spatial degrees of
freedom and is block diagonal with respect to elements.

Since forming and inverting the full coupled Boltzmann system is
impractical at scale, the standard approach is source iteration, which
updates the scattering source using the previous iterate. Starting
from $\vec{\Phi}_h^{(0)}$ (we take $\vec{0}$ so that the first sweep
gives the uncollided flux), for $n\ge 0$ we compute
\begin{equation}
    \label{eqn:si}
    \begin{split}
        \vec{A}_k \vec{\Psi}_{k,h}^{(n+1)} &= \vec{S} \vec{\Phi}_h^{(n)} + \vec{L}_k,\qquad k=1,\dots,N_Q,\\
        \vec{\Phi}_h^{(n+1)} &= \sum_{k=1}^{N_Q} w_k \vec{\Psi}_{k,h}^{(n+1)}.
    \end{split}
\end{equation}
Under the standard absorption condition
$c:= \| \sigma_s/\sigma_t \|_{L^\infty(\Omega)} <1$, this iteration is a contraction, but
its contraction factor can be arbitrarily close to one in optically
thick, highly scattering regimes, leading to slow convergence or
pseudo-convergence \cite{dsa_lifeline}. This motivates diffusion
synthetic acceleration, introduced in the next section, as an additive
correction based on a diffusion approximation of the iterative error.

In practice we apply \eqref{eqn:si} in a matrix-free manner. For each
$k$, the operator $\vec{A}_k$ can be applied and inverted via a
transport sweep based on an upwind ordering of elements. On general
unstructured or polytopic meshes this ordering can contain cycles, in
which case cycle-breaking or relaxation must be employed
\cite{VERMAAK2021109892,calloo2025cycle}. We denote the resulting
direction-dependent solve by $\texttt{sweep}_k[\cdot]$, meaning that
$\vec{X}=\texttt{sweep}_k[\vec{Y}]$ is the solution of
$\vec{A}_k\vec{X}=\vec{S}\vec{Y}+\vec{L}_k$.

\begin{remark}[Parallelisation]
For a fixed iterate $\vec{\Phi}_h^{(n)}$, the updates
$\vec{\Psi}_{k,h}^{(n+1)}$ for different ordinates are independent and
can be parallelised over $k$. Each sweep can also be parallelised in
space; see \cite{baker1998sn,calloo2025cycle} for representative
strategies.
\end{remark}

\section{Diffusion Synthetic Acceleration and its discontinuous Galerkin scheme}
\label{section:dsa_dg}
        
Source iteration can converge slowly in the diffusive regime. In the
monoenergetic, isotropically scattering setting this is characterised
by a scattering ratio $c:=\|\sigma_s/\sigma_t\|_{L^\infty(\Omega)}$
close to one together with optically thick cells, for example
$\sigma_t h_\kappa\gg 1$ on a significant portion of the mesh. In this
regime the convergence of source iteration deteriorates towards its
scattering-dominated limit, and the remaining error is dominated by
low-frequency, nearly isotropic components \cite{dsa_lifeline}.

We therefore employ diffusion synthetic acceleration (DSA), which
applies a diffusion approximation to an equation for the iterative
error and uses the resulting scalar correction as an additive
update. We first describe the diffusion approximation on the continuum
and then present its DG discretisation.

\subsection{The Diffusion Approximation}

On the continuum, source iteration for \eqref{eqn:bte} may be written
in predictor--corrector form. Given a scalar flux iterate
$\phi^{(n)}$, the transport predictor $\psi^{(n+\frac12)}$ satisfies
\begin{equation}
    \label{eqn:continuum_si}
    \begin{split}
        \vec{\omega}\cdot\nabla\psi^{(n+\frac12)}(\vec{x}, \vec{\omega}) + \sigma_t(\vec{x})\psi^{(n+\frac12)}(\vec{x}, \vec{\omega})
        &= \sigma_s(\vec{x}) \phi^{(n)}(\vec{x}) + f(\vec{x},\vec{\omega})
        \qquad \text{in }\Omega\times\mathbb{S}^{d-1},\\
        \psi^{(n+\frac12)} &= 0 \qquad \text{on }\Gamma_- .
    \end{split}
\end{equation}
We then define the intermediate scalar flux by
\[
\phi^{(n+\frac12)}(\vec{x}) := \langle \psi^{(n+\frac12)}(\vec{x},\cdot)\rangle.
\]

To derive the DSA correction, define the predictor error
\[
\varepsilon^{(n+\frac12)} := \psi - \psi^{(n+\frac12)},
\]
and the scalar error
\[
e^{(n)} := \phi - \phi^{(n)}.
\]
Subtracting \eqref{eqn:continuum_si} from the exact transport equation
yields
\begin{equation}
    \label{eqn:continuum_error_transport}
    \begin{split}
        \vec{\omega}\cdot\nabla \varepsilon^{(n+\frac12)} + \sigma_t \varepsilon^{(n+\frac12)}
        &= \sigma_s e^{(n)} \quad \text{in }\Omega\times\mathbb{S}^{d-1},\\
        \varepsilon^{(n+\frac12)} &= 0 \quad \text{on }\Gamma_- .
    \end{split}
\end{equation}
Averaging over angle gives
\[
\phi - \phi^{(n+\frac12)} = \langle \varepsilon^{(n+\frac12)} \rangle.
\]
In an optically thick, highly scattering regime, the slowly decaying
error modes are well approximated by a diffusion model. In the
isotropic setting this leads to a diffusion correction
$\delta^{(n+1)}$ defined by
\begin{equation}
    \label{eqn:dsa_continuum_correction}
    -\nabla\cdot\!\left(D(\vec{x})\nabla\delta^{(n+1)}(\vec{x})\right) + \sigma_a(\vec{x})\delta^{(n+1)}(\vec{x})
    = \sigma_s(\vec{x})\Big(\phi^{(n+\frac12)}(\vec{x}) - \phi^{(n)}(\vec{x})\Big)
    \quad \text{in }\Omega,
\end{equation}
with $\sigma_a:=\sigma_t-\sigma_s$. In the diffusion limit associated
with the source-iteration error dynamics, the isotropic diffusion
tensor is
\begin{equation}
    \label{eqn:dsa_diffusion_tensor}
    D(\vec{x}) = \frac{\sigma_s(\vec{x})}{d \sigma_t(\vec{x})^2} \mathbf{I}.
\end{equation}
The accelerated scalar update is then
\begin{equation}
    \label{eqn:dsa_update}
    \phi^{(n+1)} = \phi^{(n+\frac12)} + \delta^{(n+1)}.
\end{equation}
This diffusion approximation is key to DSA. To simplify notation, we
rewrite \eqref{eqn:dsa_continuum_correction} as
\begin{equation}
    \label{eqn:dsa_equation}
    \mathcal{D} \delta^{(n+1)} = Q^{(n)},
\end{equation}
where the diffusion operator and source term are
\begin{equation}
    \label{eqn:dsa_equations_detail}
    \mathcal{D}\varphi := -\nabla\cdot(D\nabla\varphi) + \sigma_a\varphi,
    \qquad
    Q^{(n)} := \sigma_s\left(\phi^{(n+1/2)} - \phi^{(n)}\right).
\end{equation}

The diffusion problem \eqref{eqn:dsa_equation} must be equipped with
boundary conditions compatible with the transport inflow model. In the
numerical section we consider both a strong vacuum surrogate
(Dirichlet) and the asymptotically motivated Marshak (Robin)
condition. Asymptotically correct boundary modelling requires care due
to boundary layers and geometry \cite{malvagi1991initial}; we return
to this when defining the discrete DSA operator and its weak
imposition of boundary conditions.

\subsection{Interior penalty DGFEM scheme for the diffusion equation}

We now discretise the diffusion correction problem
\eqref{eqn:dsa_equation} on the same polytopic mesh $\mathcal{T}_h$
and in the same discontinuous polynomial space as used for the
transport solve. For fixed $n$, recall that the DSA correction
$\delta^{(n+\frac12)}$ is defined by
\[
\mathcal{D} \delta^{(n+\frac12)} = Q^{(n)}
\qquad\text{in }\Omega,
\]
with $\mathcal{D}\varphi =
-\nabla\cdot(D\nabla\varphi)+\sigma_a\varphi$ and
$D=1/(d\sigma_t)>0$. Since $D>0$ and $\sigma_a\ge 0$, the operator is
elliptic and the symmetric interior penalty (SIP) formulation is
symmetric and, for sufficiently large penalty, coercive in the usual
DG energy norm. On polytopic meshes, the standard consistency and
approximation results apply under the same mesh-regularity assumptions
as in \cite{cangiani_dg} and in related DG treatments of radiative
transfer \cite{RAGUSA2015195}.

Let $\Gamma_{\mathrm{int}}$ denote the union of all interior facets of
$\mathcal{T}_h$. We partition the boundary facets into disjoint
Dirichlet and Robin sets, $\Gamma_{\mathrm{D}}$ and
$\Gamma_{\mathrm{R}}$, so that
$\Gamma=\Gamma_{\mathrm{D}}\cup\Gamma_{\mathrm{R}}$ and
$\Gamma_{\mathrm{D}}\cap\Gamma_{\mathrm{R}}=\emptyset$. We consider
boundary conditions of the form
\begin{equation}
    \label{eqn:elliptic_boundary_conditions}
    \delta = g_{\mathrm{D}}\quad\text{on }\Gamma_{\mathrm{D}},
    \qquad
    \vec{n}\cdot(D\nabla\delta) + g_{\mathrm{R}} \delta = r_{\mathrm{R}}\quad\text{on }\Gamma_{\mathrm{R}},
\end{equation}
where $\vec{n}$ is the outward unit normal on $\Gamma$ and where
$g_{\mathrm{R}}\ge 0$ is a Robin coefficient. In the numerical
experiments we specialise \eqref{eqn:elliptic_boundary_conditions} to
Dirichlet and Marshak-type conditions in
Section~\ref{section:diffusion_bcs}.

We next introduce DG trace, jump and average operators. Let $\kappa_i$
and $\kappa_j$ be neighbouring elements sharing an interior facet
$F=\partial\kappa_i\cap\partial\kappa_j\subset\Gamma_{\mathrm{int}}$,
and let $\vec{n}_i$ and $\vec{n}_j$ denote the outward unit normals on
$F$ relative to $\kappa_i$ and $\kappa_j$. For a scalar function $v$
and a vector field $\vec{q}$ that are smooth within each element, we
define the averages
\begin{equation}
    \label{eqn:averages}
    \avg{v} := \frac{1}{2}\big(v_{\kappa_i}^+ + v_{\kappa_j}^+\big),
    \qquad
    \avg{\vec{q}} := \frac{1}{2}\big(\vec{q}_{\kappa_i}^+ + \vec{q}_{\kappa_j}^+\big),
\end{equation}
and the jumps
\begin{equation}
    \label{eqn:jumps}
    \jump{v} := v_{\kappa_i}^+ \vec{n}_i + v_{\kappa_j}^+ \vec{n}_j,
    \qquad
    \jump{\vec{q}} := \vec{q}_{\kappa_i}^+\cdot\vec{n}_i + \vec{q}_{\kappa_j}^+\cdot\vec{n}_j.
\end{equation}
On a boundary facet $F\subset\Gamma$ corresponding to an element
$\kappa_i$, we use the reduced definitions
\begin{equation}
    \label{eqn:boundary_jumps_averages}
    \avg{v} := v_{\kappa_i}^+,\quad
    \avg{\vec{q}} := \vec{q}_{\kappa_i}^+,\quad
    \jump{v} := v_{\kappa_i}^+ \vec{n},\quad
    \jump{\vec{q}} := \vec{q}_{\kappa_i}^+\cdot\vec{n}.
\end{equation}

We seek a DG approximation $\delta_h^{(n+\frac12)}\in\mathbb{V}_p$ satisfying
\begin{equation}
    \label{eqn:dsa_dg}
    B\big(\delta_h^{(n+\frac12)},v_h\big) = \ell^{(n)}(v_h)\qquad \forall v_h\in\mathbb{V}_p,
\end{equation}
where the SIP bilinear form is defined by
\begin{equation}
    \label{eqn:diffusion_bilinear_form}
    \begin{split}
        B(u_h,v_h)
        :=\ &\sum_{\kappa\in\mathcal{T}_h}\int_{\kappa}\big(D\nabla u_h\cdot\nabla v_h + \sigma_a u_h v_h\big) d\vec{x}\\
        &-\int_{\Gamma_{\mathrm{int}}\cup\Gamma_{\mathrm{D}}}
        \Big(\avg{D\nabla u_h}\cdot\jump{v_h} + \avg{D\nabla v_h}\cdot\jump{u_h}
        - \sigma_{\mathrm{SIP}} \jump{u_h}\cdot\jump{v_h}\Big) ds\\
        &+\int_{\Gamma_{\mathrm{R}}} g_{\mathrm{R}} u_h v_h ds,
    \end{split}
\end{equation}
and the right-hand side is
\begin{equation}
    \label{eqn:diffusion_rhs}
    \ell^{(n)}(v_h)
    := \sum_{\kappa\in\mathcal{T}_h}\int_{\kappa} Q^{(n)} v_h d\vec{x}
    - \int_{\Gamma_{\mathrm{D}}} g_{\mathrm{D}}\Big(\vec{n}\cdot(D\nabla v_h) - \sigma_{\mathrm{SIP}} v_h\Big) ds
    + \int_{\Gamma_{\mathrm{R}}} r_{\mathrm{R}} v_h ds.
\end{equation}
The penalty parameter $\sigma_{\mathrm{SIP}}$ is chosen facetwise. Its
scaling is critical for DSA effectiveness in thick regimes and on
anisotropic polytopic meshes, and we specify our choices (including
transport-informed variants) in
Section~\ref{section:penalty_parameter}.

\subsection{Boundary Conditions in the context of Diffusion Synthetic Acceleration}
\label{section:diffusion_bcs}

The diffusion correction \eqref{eqn:dsa_equation} requires boundary
conditions. For the SI error dynamics, matching the diffusion outer
expansion to a kinetic boundary layer yields a homogeneous Robin
condition of Marshak type,
\begin{equation}
    \label{eqn:marshak_kappa}
    \vec{n}\cdot (D\nabla \varphi) + \kappa \varphi = 0 \quad \text{on }\Gamma,
    \qquad
    \kappa := \frac{1}{2}\langle |\vec{\omega}\cdot\vec{n}|\rangle.
\end{equation}
With the angular average $\langle\cdot\rangle$ and the normalisation
$\sum_k w_k=1$, this constant evaluates to $\kappa=1/4$ in three
dimensions and $\kappa=1/\pi$ in two dimensions.

In the experiments below we consider two practical choices for the
diffusion correction boundary treatment.

\subsubsection{Homogeneous Dirichlet boundary conditions}

We impose
\[
\varphi = 0 \quad \text{on }\Gamma,
\]
that is, $\Gamma_{\mathrm D}=\Gamma$ and $\Gamma_{\mathrm
  R}=\emptyset$ in \eqref{eqn:elliptic_boundary_conditions}. This is a
common computational choice and is straightforward to implement in the
interior penalty formulation.

\subsubsection{Homogeneous Robin (Marshak) boundary conditions}

We impose the Marshak condition \eqref{eqn:marshak_kappa} on the full
boundary, that is, $\Gamma_{\mathrm R}=\Gamma$ and $\Gamma_{\mathrm
  D}=\emptyset$. In the notation of
\eqref{eqn:elliptic_boundary_conditions} this corresponds to
\[
g_{\mathrm R}=\kappa,
\qquad
r_{\mathrm R}=0,
\]
so that
\[
\vec{n}\cdot(D\nabla\varphi) + \kappa \varphi = 0
\quad \text{on }\Gamma.
\]

\subsection{Modified interior penalty scheme}
\label{section:penalty_parameter}

The SIP formulation \eqref{eqn:diffusion_bilinear_form} is written
with the standard (arithmetic) averages and jumps from
\eqref{eqn:averages}--\eqref{eqn:jumps}. In particular, we do not use
diffusion-balanced (weighted) averages in the consistency terms. The
stabilisation is therefore controlled entirely by the choice of the
penalty parameter on facets.

A standard polytopic SIP choice is
\begin{equation}
    \label{eqn:sip_penalty_polytopic}
    \sigma_{\mathrm{SIP}}(\vec{x}) := 
    \begin{cases}
        C\max_{\kappa\in\{\kappa^-, \kappa^+\}}\left\{A_F|_{\kappa}  C_{\mathrm{INV}}(p, \kappa, F) \dfrac{p^2|F|}{|\kappa|}\right\}
        &\quad \vec{x}\in F\in\Gamma_{\mathrm{int}},\ F\subset\partial\kappa^-\cap\partial\kappa^+, \\[1.0em]
        C A_F  C_{\mathrm{INV}}(p, \kappa, F) \dfrac{p^2|F|}{|\kappa|}
        &\quad \vec{x}\in F\in\Gamma_{\mathrm{D}},\ F\subset\partial\kappa,
    \end{cases}
\end{equation}
where $A_F:=\|\sqrt{D} \vec{n}\|_{\leb{\infty}(F)}^2$ for each facet
$F\subset\partial\kappa$,
$F\in\Gamma_{\mathrm{int}}\cup\Gamma_{\mathrm{D}}$, and $C>0$ is taken
sufficiently large. The inverse constant $C_{\mathrm{INV}}$ is given
by
\begin{equation}
    \label{eqn:cinv_polytopic}
    C_{\mathrm{INV}}(p, \kappa, F) := C_{\mathrm{inv}}\min\left\{\frac{|\kappa|}{\Lambda_F(\kappa)},  p^{2(d-1)}\right\},
\end{equation}
under the usual $p$-coverability assumptions; see
\cite[Lemma~11]{cangiani_dg}.

In the thick diffusive regime relevant to DSA, the diffusion
coefficient $D$ is small (and $\sigma_a$ may be small), so
$\sigma_{\mathrm{SIP}}$ can become too weak to control jump modes in a
manner comparable with the upwind dissipation present in the transport
discretisation. The MIP remedy is to enforce a transport-scale floor
on the penalty, face by face, using the same $S_N$ quadrature used in
the transport solve. We assume the quadrature weights are normalised
so that $\sum_{m=1}^{N_Q} w_m = 1$.

For an interior facet $F\in\Gamma_{\mathrm{int}}$, fix any unit normal
$\vec{n}_F$ (its orientation is irrelevant below) and define the
transport face moment
\begin{equation}
    \label{eqn:transport_face_moment_interior}
    C_{0,F} := \frac{1}{2}\sum_{m=1}^{N_Q} w_m \big|\vec{\omega}_m\cdot \vec{n}_F\big|.
\end{equation}
For a boundary facet $F\subset\Gamma$, define the outflow moment
\begin{equation}
    \label{eqn:transport_face_moment_boundary}
    C_{0,F}^{\partial} := \sum_{m=1}^{N_Q} w_m \max\{0,\vec{\omega}_m\cdot \vec{n}\}.
\end{equation}
We then define the modified interior penalty as
\begin{equation}
    \label{eqn:mip_penalty}
    \sigma_{\mathrm{MIP}}(\vec{x}) :=
    \begin{cases}
        \max\{\sigma_{\mathrm{SIP}}(\vec{x}),  C_{0,F}\}, & \vec{x}\in F\in\Gamma_{\mathrm{int}},\\
        \max\{\sigma_{\mathrm{SIP}}(\vec{x}),  C_{0,F}^{\partial}\}, & \vec{x}\in F\in\Gamma_{\mathrm{D}}.
    \end{cases}
\end{equation}
In \eqref{eqn:diffusion_bilinear_form} we replace
$\sigma_{\mathrm{SIP}}$ by $\sigma_{\mathrm{MIP}}$. This preserves the
standard SIP structure (with arithmetic averages) while preventing the
penalty from degenerating relative to the upwind transport dissipation
in optically thick, highly scattering settings.

\begin{remark}
The constants in
\eqref{eqn:transport_face_moment_interior}--\eqref{eqn:transport_face_moment_boundary}
are the discrete analogues of the half-range angular moments that
appear in transport energy identities.  They depend only on the chosen
quadrature and face orientation, not on $D$.
\end{remark}
            
\subsection{Implementation and matrix notation}

We now summarise the DSA-accelerated iteration in a form suitable for
implementation. As on the continuum, the method consists of two steps, a transport sweep (source iteration) producing an intermediate scalar
flux, followed by a diffusion solve producing an additive correction.

Let $\vec{\Psi}_{k,h}$ denote the coefficient vector of the discrete
ordinate solution $\psi_{k,h}\in\mathbb{V}_p$ and let $\vec{\Phi}_h$
denote the coefficient vector of the scalar flux
$\phi_h=\sum_{k=1}^{N_Q}w_k\psi_{k,h}$. As in \eqref{eqn:si_full}, the
discrete transport systems are
\[
\vec{A}_k\vec{\Psi}_{k,h}=\vec{S}\vec{\Phi}_h+\vec{L}_k,\qquad k=1,\dots,N_Q,
\]
where $\vec{A}_k$ is the upwind DG transport matrix, $\vec{S}$ is the
(spatial) reaction matrix associated with $\sigma_s$ (so that
$\vec{S}\vec{\Phi}_h$ represents $\int_\Omega \sigma_s \phi_h v_h$ in
coefficient form), and $\vec{L}_k$ is the load vector.

Given $\vec{\Phi}_h^{(n)}$, a transport sweep produces directional
iterates $\vec{\Psi}_{k,h}^{(n+1)}$ and hence the intermediate scalar
flux
\[
\vec{\Phi}_h^{(n+\frac12)} := \sum_{k=1}^{N_Q} w_k \vec{\Psi}_{k,h}^{(n+1)}.
\]
The DSA correction is defined by the discrete diffusion solve: find
$\vec{\delta}_h^{(n+1)}$ such that
\begin{equation}
    \label{eqn:dsa_matrix}
    \vec{D} \vec{\delta}_h^{(n+1)} = \vec{S}\big(\vec{\Phi}_h^{(n+\frac12)}-\vec{\Phi}_h^{(n)}\big),
\end{equation}
where $\vec{D}$ is the matrix corresponding to the SIP/MIP
discretisation \eqref{eqn:dsa_dg}--\eqref{eqn:diffusion_bilinear_form}
(with the chosen boundary treatment from
Section~\ref{section:diffusion_bcs} and penalty from
Section~\ref{section:penalty_parameter}). The accelerated scalar
update is then
\begin{equation}
    \label{eqn:dsa_matrix_update}
    \vec{\Phi}_h^{(n+1)} = \vec{\Phi}_h^{(n+\frac12)} + \vec{\delta}_h^{(n+1)}.
\end{equation}

We present the resulting algorithm in Algorithm~\ref{alg:dsa}. The
diffusion matrix $\vec{D}$ is independent of the iteration index and
is assembled once. Only the right-hand side in \eqref{eqn:dsa_matrix}
must be formed at each iteration, which is inexpensive compared with
transport sweeps.

\begin{algorithm}
    \caption{DSA iteration scheme}\label{alg:dsa}
    \KwData{Coefficients $\sigma_t,\sigma_s$, source $f$ and inflow data $g_{\mathrm D}$; angular quadrature $\{(w_k,\vec{\omega}_k)\}_{k=1}^{N_Q}$; tolerance $\epsilon$}
    $\vec{D} \gets \texttt{assemble\_diffusion\_matrix}(\sigma_t,\sigma_s)$\Comment*[r]{SIP/MIP diffusion operator}
    
    $\vec{\Phi}_h^{(0)} \gets \vec{0}$\Comment*[r]{Initial scalar flux}
    $n \gets 0$, $r \gets \infty$
    
    \While{$r>\epsilon$}{
        \For{$k \gets 1$ to $N_Q$}{
            $\vec{\Psi}_{k,h}^{(n+\frac12)} \gets \texttt{sweep}_k[\vec{\Phi}_h^{(n)}]$\Comment*[r]{Transport predictor}
        }
        $\vec{\Phi}_h^{(n+\frac12)} \gets \sum_{k=1}^{N_Q} w_k \vec{\Psi}_{k,h}^{(n+\frac12)}$\Comment*[r]{Intermediate scalar flux}
        
        $\vec{b}^{(n)} \gets \vec{S}\big(\vec{\Phi}_h^{(n+\frac12)}-\vec{\Phi}_h^{(n)}\big)$\Comment*[r]{DSA source}
        $\vec{\delta}_h^{(n+1)} \gets \texttt{solve\_diffusion}(\vec{D},\vec{b}^{(n)})$\Comment*[r]{Diffusion correction}
        
        $\vec{\Phi}_h^{(n+1)} \gets \vec{\Phi}_h^{(n+\frac12)} + \vec{\delta}_h^{(n+1)}$\Comment*[r]{Accelerated update}
        
        $r \gets \|\vec{\Phi}_h^{(n+1)}-\vec{\Phi}_h^{(n)}\| / \|\vec{\Phi}_h^{(n+1)}\|$
        $n \gets n + 1$
    }
    \Return $\{\vec{\Psi}_{k,h}^{(n-\frac12)}\}_{k=1}^{N_Q}$ and $\vec{\Phi}_h^{(n)}$
\end{algorithm}

The matrix $\vec{D}$ is symmetric positive definite for both SIP and
MIP choices, so conjugate gradient methods (optionally preconditioned,
for example by AMG) are natural candidates for
\texttt{solve\_diffusion}. In the numerical experiments of this paper
we use a sparse PLU factorisation for the diffusion solve in order to
obtain consistent and directly comparable timings across different
parameter choices; replacing this with CG/AMG changes the absolute
wall-clock times but does not affect the measured contraction factors
of the outer DSA iteration.

\section{Numerical Experiments}
\label{section:numerics}

We now present numerical results for the DSA schemes introduced in \S
\ref{section:dsa_dg}. The aim is to quantify the effectiveness of DSA
in adverse diffusive regimes and to assess its behaviour under
standard advanced DGFEM practices, in particular
$hp$-refinement. Since DSA performance is strongly influenced by the
diffusion interior penalty, we focus on parameters that directly
affect this term, notably cell sizes and local polynomial degree.

To compare the iterative schemes, we report an \emph{empirical
convergence factor} for the outer iteration. Let $\phi^{(n)}$ denote
the scalar flux iterate at outer iteration $n$ and let $\phi^\ast$
denote a high-accuracy discrete reference solution of the fully
discretised problem. We define the per-iterate error reduction factor
by
\begin{equation}
    \label{eqn:rho_n}
    \rho^{(n)} = \frac{\|\phi^{(n+1)} - \phi^\ast\|_{L^2(\Omega)}}{\|\phi^{(n)} - \phi^\ast\|_{L^2(\Omega)}},
\end{equation}
and report the geometric mean
\begin{equation}
    \label{eqn:rho_geom}
    \rho = \left(\prod_{n=0}^{M-1}\rho^{(n)}\right)^{\frac{1}{M}}.
\end{equation}
Here $M$ is chosen so that $\phi^{(M)}$ is the last iterate whose
error remains above $10^{-10}$ in $L^2(\Omega)$, that is,
$\|\phi^{(M+1)}-\phi^\ast\|_{L^2(\Omega)}<10^{-10}$. In this way we
exclude the final iterations approaching machine precision while still
capturing the observed convergence behaviour over the main part of the
iteration. If $\rho>1$, we classify the scheme as divergent.

Unless otherwise stated, all experiments use the same baseline
configuration. We solve both the $S_N$ transport problem and the DSA
diffusion problem on the domain $\Omega=(0,10)^2$ with a bounded
Voronoi tessellation containing 1024 elements. The mesh is smoothed
using Lloyd's algorithm to reduce element anisotropy and to isolate
the effect of other parameters. The resulting tessellations typically
have around six facets per element, with a minimum of four and a
maximum of eight. The spatial approximation uses the same
discontinuous polynomial space for transport and diffusion with
polynomial degree $p=1$ by default.

We take constant cross-sections $\sigma_t$ and $\sigma_s$ (homogeneous
material). To explore optical thickness, we vary $\sigma_t$ and set
$\sigma_s=c\sigma_t$ so that $\sigma_a=(1-c)\sigma_t$. Unless
specified, we use scattering ratio $c=0.999$. Angular discretisation
uses $N_Q=16$ discrete ordinates based on the trapezoidal rule on
$\mathbb{S}^1$: for $m=1,\ldots,N_Q$, $\theta_m=2\pi(m-1)/N_Q$,
$\vec{\omega}_m=(\cos\theta_m,\sin\theta_m)$ and
$w_m=1/N_Q$, so that $\sum_m w_m=1$ and
$\sum_m w_m\vec{\omega}_m=\vec{0}$.

To eliminate the influence of inner iteration errors from the
diffusion correction, we solve the DSA linear system directly using a
sparse LU factorisation from \texttt{scipy}. By contrast, the
transport solves are performed by directional sweeps, as described in
\S\ref{section:transport_dg}. This separation ensures that the
reported convergence factors reflect the outer SI/DSA iteration
itself, rather than the stopping tolerances or preconditioning choices
of an inner iterative diffusion solver.  Replacing the direct
diffusion solve by CG, AMG or related iterative methods would change
absolute wall-clock times, but not the observed outer-iteration
convergence factors provided the diffusion problem is solved to
sufficient accuracy. The outer SI/DSA iteration is run until the
relative scalar flux change satisfies $r^{(n)}<10^{-12}$ or until 50
iterations are reached. We initialise with
$\vec{\Phi}_h^{(0)}=\vec{0}$, so that the first transport sweep
produces the uncollided flux.

The SIP penalty contains mesh- and degree-dependent constants that are
not available in closed form for general polytopal elements. We
therefore absorb these constants into a single prefactor and set this
prefactor equal to $10$, which is a standard practical choice in SIP
discretisations for elliptic problems. More robust alternatives exist;
see \cite{dong2022robust} and the references therein.

To enable consistent comparisons across parameter variations, we
employ a manufactured solution in two dimensions,
\begin{equation}
    \label{eqn:example_solution}
    \psi(\vec{x}, \vec{\omega}) = \omega_{x_1}^2\sin{\pi x_1}\sin{\pi x_2},
    \qquad \vec{x}=(x_1,x_2)\in(0,10)^2,
\end{equation}
with $\vec{\omega}=(\omega_{x_1},\omega_{x_2})$. The corresponding
ordinate-wise solutions are
\begin{equation}
    \label{eqn:example_solution_ord}
    \psi_k(\vec{x}) = \omega_{k,x_1}^2\sin{\pi x_1}\sin{\pi x_2},
    \qquad \vec{x}\in(0,10)^2,
\end{equation}
where $\vec{\omega}_k=(\omega_{k,x_1},\omega_{k,x_2})$, and the source
term and inflow data are chosen consistently so that
\eqref{eqn:example_solution} satisfies the transport model exactly.
This choice provides a non-trivial spatial and angular dependence
while remaining smooth enough to support systematic studies of angular
refinement, mesh refinement and polynomial degree.

Across the numerical section we test four DSA variants, obtained by
combining SIP or MIP diffusion discretisations with either homogeneous
Dirichlet or Marshak (Robin) diffusion boundary conditions. We also
report results for unaccelerated source iteration as a baseline.

\subsection{Experiment 1: Comparison of boundary conditions and interior penalty schemes}
\label{subsec:exp1_bc_penalty}

In this experiment we study how the total macroscopic cross-section
$\sigma_t$ influences the convergence of the outer iteration, and how
this depends on the diffusion boundary treatment and the choice of
interior penalty scheme. We compare the unaccelerated source iteration
baseline with four DSA variants, SIP--Dirichlet, SIP--Marshak (Robin),
MIP--Dirichlet and MIP--Marshak.

\begin{figure}[htbp]
    \centering
    \begin{tikzpicture}
    \begin{axis}[
        width=0.9\textwidth,
        height=0.5\textwidth,
        ymin = 0,
        ymax = 1,
        xmax = 1e6,
        very thick,
        xmode=log,
        xlabel= {Total macroscopic cross-section, $\sigma_t$},
        ylabel= {Empirical convergence factor, $\rho$},
        grid=both,
        minor grid style={gray!25},
        major grid style={gray!25},
        legend pos=north west,
      ]
      \addlegendimage{empty legend}
      \addlegendentry{\hspace{-.6cm}\textbf{Acceleration}}
      \addlegendimage{empty legend}
      \addlegendentry{\hspace{-.8cm}\textbf{method}}

      \addplot[color={rgb,255:red,0;green,114;blue,178}] table[x=scalar, y=spectral_radius, col sep=comma] {tikz_Schaubilder/sr_calculations/_new_results/None.csv};
      \addlegendentry{None};
      \addplot[color={rgb,255:red,213;green,94;blue,0}, dashed] table[x=scalar, y=spectral_radius, col sep=comma] {tikz_Schaubilder/sr_calculations/_new_results/Dirichlet_DSA.csv};
      \addlegendentry{Dirichlet SIP};
      \addplot[color={rgb,255:red,0;green,158;blue,115}, dashed] table[x=scalar, y=spectral_radius, col sep=comma] {tikz_Schaubilder/sr_calculations/_new_results/Robin_DSA.csv};
      \addlegendentry{Robin SIP};
      \addplot[color={rgb,255:red,230;green,159;blue,0}, dash pattern=on 4pt off 2pt on 1pt off 2pt] table[x=scalar, y=spectral_radius, col sep=comma] {tikz_Schaubilder/sr_calculations/_new_results/Dirichlet_MIP_DSA.csv};
      \addlegendentry{Dirichlet MIP};
      \addplot[color={rgb,255:red,86;green,180;blue,233}, dash pattern=on 4pt off 2pt on 1pt off 2pt] table[x=scalar, y=spectral_radius, col sep=comma] {tikz_Schaubilder/sr_calculations/_new_results/Robin_MIP_DSA.csv};
      \addlegendentry{Robin MIP};
  
    \end{axis}
  \end{tikzpicture}
    \caption{Experiment~\ref{subsec:exp1_bc_penalty}. Empirical
      convergence factor as a function of the total macroscopic
      cross-section $\sigma_t$, comparing unaccelerated source
      iteration with DSA using SIP or MIP diffusion discretisations
      and either homogeneous Dirichlet or Marshak (Robin) diffusion
      boundary conditions.}
    \label{fig:exp1_comparison}
\end{figure}
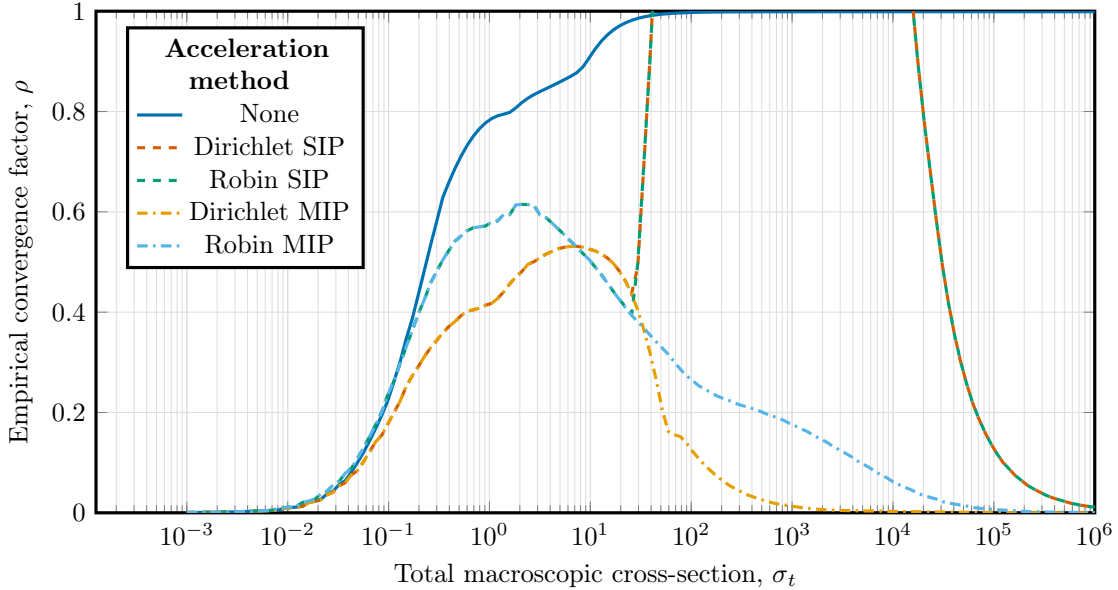

We use the manufactured solution \eqref{eqn:example_solution} and
solve both the transport and diffusion correction problems on the same
mesh and in the same finite element space. The resulting observed
convergence factors, computed as described in
Section~\ref{section:numerics}, are reported in
Figure~\ref{fig:exp1_comparison}.

Figure~\ref{fig:exp1_comparison} shows as $\sigma_t$ increases and the
regime becomes optically thick and highly scattering, convergence
deteriorates and the observed convergence factor approaches one. All
DSA variants substantially improve convergence in the intermediate and
diffusive regimes. However, the SIP-based DSA schemes exhibit a clear
loss of robustness. The convergence improves up to an intermediate
regime (here around $\sigma_t\approx 17$), after which the SIP
variants lose stability and eventually diverge. In contrast, the
MIP-based DSA schemes remain convergent for all values of $\sigma_t$
tested, consistent with the role of MIP in preventing
under-penalisation in optically thick settings.

The diffusion boundary conditions also affect performance, most
clearly in the intermediate regime. For smaller $\sigma_t$, the
homogeneous Dirichlet conditions yields improved convergence relative to
Marshak (Robin) conditions. Once the SIP variants fail, the two
MIP variants accelerate convergence more effectively, and for large $\sigma_t$
their observed convergence factors approach one another. We conjecture
that in the diffusive limit the boundary layer becomes confined to a
shrinking region near $\Gamma$, so that the influence of the precise
diffusion boundary model on the bulk convergence diminishes
\cite{malvagi1991initial}.

\subsection{Experiment 2: Comparison of the computational costs of SI and DSA}
\label{subsec:exp2_costs}

Having established that DSA can substantially reduce the observed
convergence factor in diffusive regimes
(Experiment~\ref{subsec:exp1_bc_penalty}), we now examine the
computational cost of including a diffusion correction at each outer
iteration. Although the diffusion matrix $\vec{D}$ is assembled once,
each outer iteration requires formation of a new diffusion source term
and the solution of an auxiliary linear system. This subsection
quantifies the resulting wall-clock overhead and its dependence on the
number of discrete ordinates.

We focus on the \emph{online} costs during the iterative solve. We
measure (i) the average per-iteration transport time (the sweep step),
(ii) the average per-iteration diffusion cost, comprising diffusion
source assembly and diffusion solve, and (iii) the total
per-iteration time, which additionally includes vector updates and
reductions.

For clarity we report results for a single representative accelerated
scheme, namely MIP with homogeneous Dirichlet diffusion boundary
conditions. We also perform the transport solves sequentially, so that
the reported transport--diffusion cost balance is not dominated by
hardware-specific parallelisation effects.

The key observation is that DSA can be relatively expensive per
iteration when very small quadrature sets are used. For sufficiently
small $N_Q$, the diffusion correction dominates the per-iteration
time. As $N_Q$ increases, the transport work scales approximately with
the number of ordinates, whereas the diffusion work remains
essentially constant, so the transport sweeps increasingly dominate
the per-iteration cost and the DSA overhead is amortised. This is
particularly relevant in settings with non-trivial angular structure,
where larger quadrature sets are required for accuracy and the reduced
iteration count from DSA can yield a clear overall benefit.

\begin{figure}[htbp]
    \centering
    \begin{tikzpicture}
    \begin{groupplot}[
        group style={
            group size=2 by 1,
            horizontal sep= 2.5cm,
            vertical sep=2.5cm,
            every plot/.style={
            very thick,
            legend pos=north west,
            },
        },
        enlargelimits=false,
    ]
 
    % Top plot
    \nextgroupplot[
        width=0.46\textwidth,
        height=0.4\textwidth,
        ymax=4, 
        xmin=2, xmax=2048,
        xmode=log, ymode=log,
        xlabel={Number of discrete ordinates, $N_Q$},
        ylabel={Wall-clock time, [s]},
        grid=major,
        title={Wall-clock time per accelerated iteration},
    ]
    
    \addplot[brown, mark options={black, scale=0.75}] 
    table[x=ndo, y=iteration_times, col sep=comma] {tikz_Schaubilder/src_t/_new_results/Dirichlet_MIP_0.5.csv};
    \addlegendentry{Total};

    \addplot[blue] table[x=ndo, y=si_times, col sep=comma] {tikz_Schaubilder/src_t/_new_results/Dirichlet_MIP_0.5.csv};
    \addlegendentry{SI};

    \addplot[green] table[x=ndo, y=dsa_times, col sep=comma] {tikz_Schaubilder/src_t/_new_results/Dirichlet_MIP_0.5.csv};
    \addlegendentry{DSA};

    \nextgroupplot[
        width=0.46\textwidth,
        height=0.4\textwidth,
        ymin=0,ymax=100, 
        xmin=2, xmax=2048,
        xmode=log,
        grid=major,
        xlabel={Number of discrete ordinates, $N_Q$},
        ylabel={Time, [\%]},
        title={Percentage time per accelerated iteration},
        legend style={at={(0.65,0.65)}}
    ]

    \addplot[blue] table[x=ndo, y expr=100*\thisrow{si_times}/\thisrow{iteration_times}, col sep=comma] {tikz_Schaubilder/src_t/_new_results/Dirichlet_MIP_0.5.csv};
    \addlegendentry{SI};
    
    \addplot[green] table[x=ndo, y expr=100*\thisrow{dsa_times}/\thisrow{iteration_times}, col sep=comma] {tikz_Schaubilder/src_t/_new_results/Dirichlet_MIP_0.5.csv};
    \addlegendentry{DSA};

    \end{groupplot}
    
\end{tikzpicture}
    \caption{Experiment~\ref{subsec:exp2_costs}. Per-iteration cost of
      DSA (MIP--Dirichlet) as a function of the number of discrete
      ordinates $N_Q$. Shown are the mean wall-clock time spent in the
      diffusion correction (source assembly plus diffusion solve) and
      the corresponding percentage of the total per-iteration time.}
    \label{fig:exp2_t_iteration}
\end{figure}
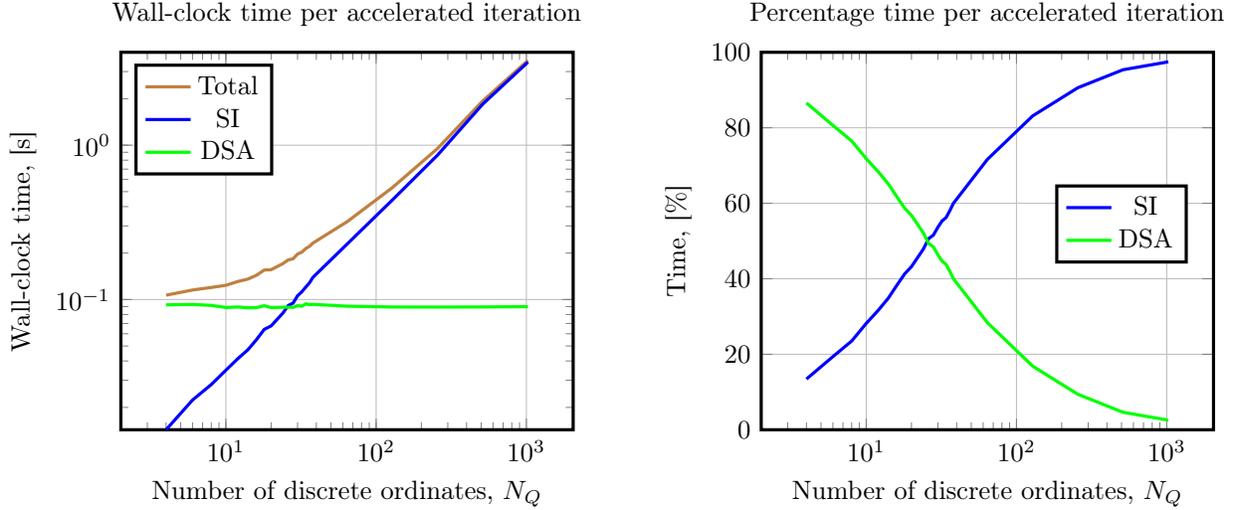

We first report the per-iteration overhead.
Figure~\ref{fig:exp2_t_iteration} shows both the average wall-clock
time attributable to the DSA correction and the corresponding
percentage of the total per-iteration time, as a function of $N_Q$.

The balance point in Figure~\ref{fig:exp2_t_iteration}, at which the
transport sweep and diffusion correction contribute comparable
fractions of the iteration time, occurs at $N_Q=26$ in this
example. At this point the diffusion correction accounts for $49\%$ of
the per-iteration time, while at $N_Q=128$ this drops to $17\%$.

\begin{figure}[htbp]
    \centering
    \begin{tikzpicture}
    \begin{groupplot}[
        group style={
            group size=2 by 1,
            horizontal sep= 2.5cm,
            vertical sep=2.5cm,
            every plot/.style={
            very thick,
            legend pos=north west,
            },
        },
        enlargelimits=false,
    ]

    \nextgroupplot[
        width=0.46\textwidth,
        height=0.4\textwidth,
        ymin=1, ymax=330,
        xmin=2, xmax=2048,
        xmode=log, ymode=log,
        grid=major,
        xlabel={Number of discrete ordinates, $N_Q$},
        ylabel={Wall-clock time, [s]},
        title={Total wall-clock time},
    ]
    
    \addplot[blue, mark options={black, scale=0.75}] 
    table[x=ndo, y=total_times, col sep=comma] {tikz_Schaubilder/src_t/_new_results/SI_0.5.csv};
    \addlegendentry{SI};

    \addplot[green, mark options={black, scale=0.75}] 
    table[x=ndo, y=total_times, col sep=comma] {tikz_Schaubilder/src_t/_new_results/Dirichlet_MIP_0.5.csv};
    \addlegendentry{DSA};

    \nextgroupplot[
        width=0.46\textwidth,
        height=0.4\textwidth,
        xmin=2, xmax=2048,
        ymin=0, ymax=330,
        xmode=log,
        grid=major,
        xlabel={Number of discrete ordinates, $N_Q$},
        ylabel={Speed-up, [\%]},
        title={Total speed-up},
    ]
    
    \addplot[violet] table[x=ndo, y=speedup, col sep=comma] {tikz_Schaubilder/src_t/_new_results/speedup_0.5.csv};

    \end{groupplot}
    
\end{tikzpicture}
    \caption{Experiment~\ref{subsec:exp2_costs}. Total wall-clock time
      to convergence for unaccelerated source iteration (SI) and DSA
      (MIP--Dirichlet), shown as a function of the number of discrete
      ordinates $N_Q$. Also shown is the percentage speed-up achieved
      by DSA.}
    \label{fig:exp2_t_total}
\end{figure}
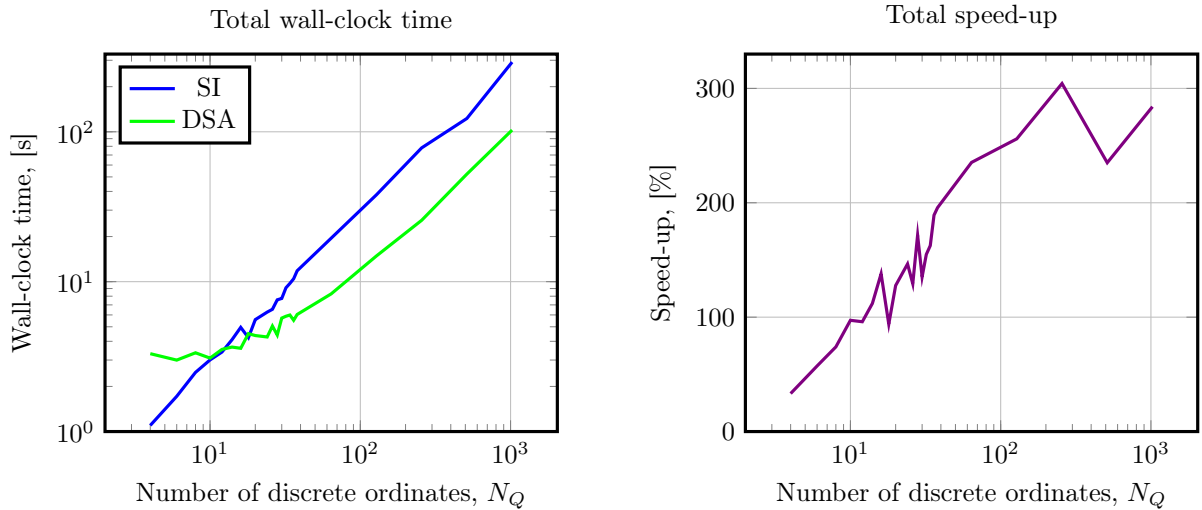

We next consider the total time to convergence.
Figure~\ref{fig:exp2_t_total} reports the total wall-clock time
required for SI and for DSA (again using MIP--Dirichlet), together
with the percentage speed-up achieved by DSA, defined by
\[
100\times\frac{T_{\mathrm{SI}}-T_{\mathrm{DSA}}}{T_{\mathrm{SI}}}.
\]

Figure~\ref{fig:exp2_t_total} shows that DSA becomes advantageous in
total runtime once the reduction in iteration count compensates for
the additional per-iteration diffusion overhead. In this example, the
SI and DSA total times coincide at $N_Q=14$, after which the speed-up
increases steadily with $N_Q$ (this is when $\sigma_t=0.5$, see Figure~\ref{fig:exp1_comparison} for observed convergence factors).

Finally, we note that the near-independence of the DSA cost with
respect to $N_Q$ assumes a fixed diffusion mesh and a fixed diffusion
solver configuration. In heterogeneous settings where $\sigma_t$
varies strongly in space, the efficiency of the diffusion solve (and
therefore the overall DSA overhead) can change. We observed this
effect in additional tests and refer to \cite{Southworth01022021} for
a detailed discussion in heterogeneous media.

\subsection{Experiment 3: Varying the scattering ratio}
\label{subsec:exp3_scattering_ratio}

A natural follow-on study is to examine how the accelerated schemes
behave as the scattering ratio approaches unity. In classical source
iteration theory, the convergence rate deteriorates as the scattering
ratio increases, but it is not a priori clear how this dependence
carries over to DSA in the diffusive limit, where $\sigma_a$ is small
and $\sigma_t$ is large.

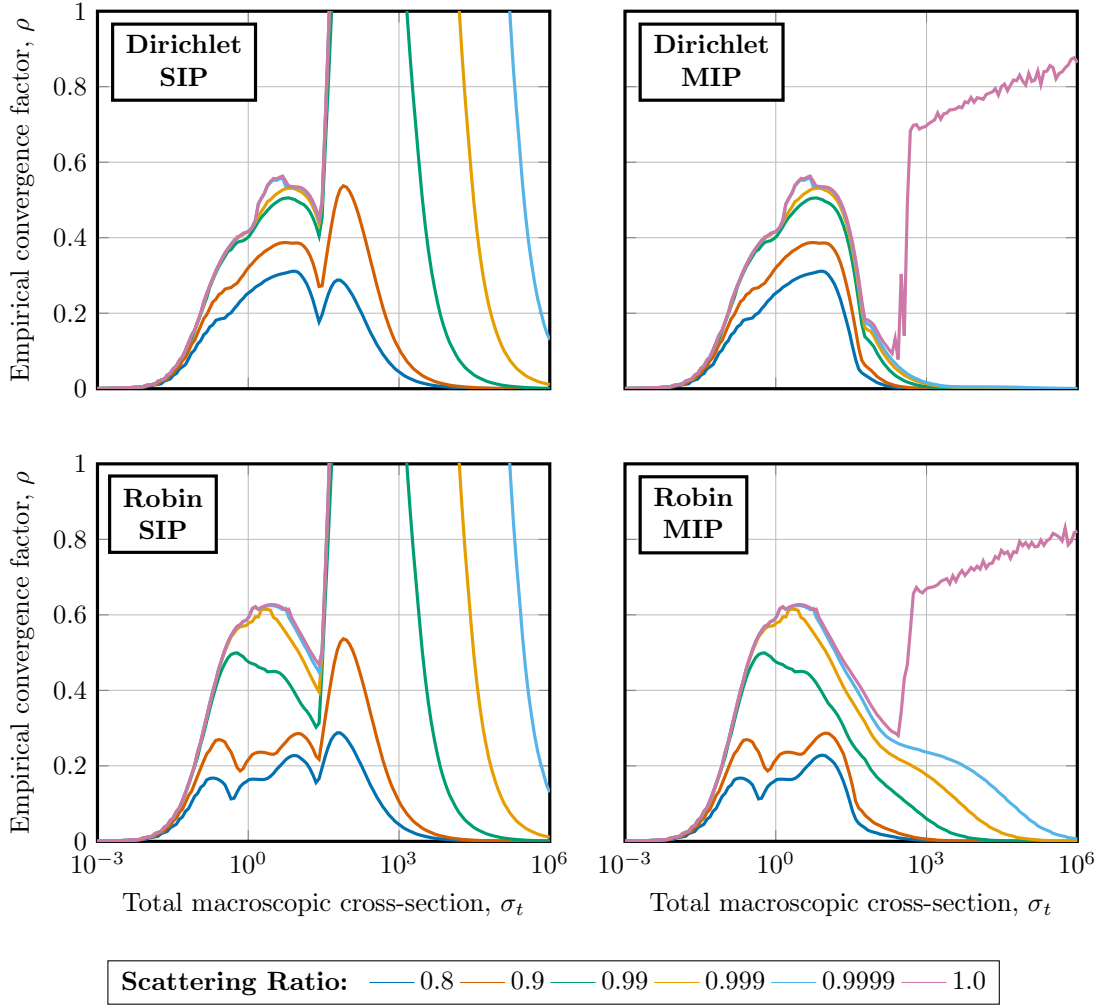
\begin{figure}[htbp]
    \centering
    \begin{tikzpicture}
    % \begin{groupplot}[
    %     group style={
    %         group size=2 by 2,
    %         horizontal sep= 1.0cm,
    %         vertical sep=1cm,
    %         x descriptions at=edge bottom,
    %         y descriptions at=edge left,
    %         every plot/.style={
    %             very thick,
    %             ymin=0, ymax=1, 
    %             xmin=1e-3, xmax=1e6,
    %             xtick={1e-3, 1e0, 1e3, 1e6},
    %             legend pos=north west,
    %         },
    %     },
    %     enlargelimits=false,
    %     xlabel={Total macrosopic cross-section, $\sigma_t$},
    %     ylabel={Spectral radius, $\rho$},
    % ]

    \begin{groupplot}[
        group style={
            group size=2 by 2,
            horizontal sep=1.0cm,
            vertical sep=1cm,
            x descriptions at=edge bottom,
            y descriptions at=edge left,
        },
        every axis/.style={
            very thick,
            grid=major,
            xmode=log,
            ymin=0, ymax=1,
            xmin=1e-3, xmax=1e6,
            xtick={1e-3,1e0,1e3,1e6},
            enlargelimits=false,
            unbounded coords=discard,
            filter discard warning=false,
            restrict y to domain=0:100,
            legend pos=north west,
        },
        width=0.46\textwidth,
        height=0.4\textwidth,
        xlabel={Total macroscopic cross-section, $\sigma_t$},
        ylabel={Empirical convergence factor, $\rho$},
    ]
 
    % Top plot
    \nextgroupplot
    
    \addlegendimage{empty legend}
    \addlegendentry{\textbf{Dirichlet}}
    \addlegendimage{empty legend}
    \addlegendentry{\textbf{SIP}}
    
    \addplot[color={rgb,255:red,0;green,114;blue,178}, mark options={black, scale=0.75}] 
    table[x=scalar, y=spectral_radius, col sep=comma] {tikz_Schaubilder/src_c/_new_results/Dirichlet_DSA_0.8.csv};

    \addplot[color={rgb,255:red,213;green,94;blue,0}, mark options={black, scale=0.75}] 
    table[x=scalar, y=spectral_radius, col sep=comma] {tikz_Schaubilder/src_c/_new_results/Dirichlet_DSA_0.9.csv};
    
    \addplot[color={rgb,255:red,0;green,158;blue,115}, mark options={black, scale=0.75}] 
    table[x=scalar, y=spectral_radius, col sep=comma] {tikz_Schaubilder/src_c/_new_results/Dirichlet_DSA_0.99.csv};
    
    \addplot[color={rgb,255:red,230;green,159;blue,0}, mark options={black, scale=0.75}] 
    table[x=scalar, y=spectral_radius, col sep=comma] {tikz_Schaubilder/src_c/_new_results/Dirichlet_DSA_0.999.csv};
    
    \addplot[color={rgb,255:red,86;green,180;blue,233}, mark options={black, scale=0.75}] 
    table[x=scalar, y=spectral_radius, col sep=comma] {tikz_Schaubilder/src_c/_new_results/Dirichlet_DSA_0.9999.csv};

    \addplot[color={rgb,255:red,204;green,121;blue,167}, mark options={black, scale=0.75}] 
    table[x=scalar, y=spectral_radius, col sep=comma] {tikz_Schaubilder/src_c/_new_results/Dirichlet_DSA_1.0.csv};
    
    \nextgroupplot
    
    \addlegendimage{empty legend}
    \addlegendentry{\textbf{Dirichlet}}
    \addlegendimage{empty legend}
    \addlegendentry{\textbf{MIP}}
    
    \addplot[color={rgb,255:red,0;green,114;blue,178}, mark options={black, scale=0.75}] 
    table[x=scalar, y=spectral_radius, col sep=comma] {tikz_Schaubilder/src_c/_new_results/Dirichlet_MIP_DSA_0.8.csv};

    \addplot[color={rgb,255:red,213;green,94;blue,0}, mark options={black, scale=0.75}] 
    table[x=scalar, y=spectral_radius, col sep=comma] {tikz_Schaubilder/src_c/_new_results/Dirichlet_MIP_DSA_0.9.csv};
    
    \addplot[color={rgb,255:red,0;green,158;blue,115}, mark options={black, scale=0.75}] 
    table[x=scalar, y=spectral_radius, col sep=comma] {tikz_Schaubilder/src_c/_new_results/Dirichlet_MIP_DSA_0.99.csv};
    
    \addplot[color={rgb,255:red,230;green,159;blue,0}, mark options={black, scale=0.75}] 
    table[x=scalar, y=spectral_radius, col sep=comma] {tikz_Schaubilder/src_c/_new_results/Dirichlet_MIP_DSA_0.999.csv};
    
    \addplot[color={rgb,255:red,86;green,180;blue,233}, mark options={black, scale=0.75}] 
    table[x=scalar, y=spectral_radius, col sep=comma] {tikz_Schaubilder/src_c/_new_results/Dirichlet_MIP_DSA_0.9999.csv};

    \addplot[color={rgb,255:red,204;green,121;blue,167}, mark options={black, scale=0.75}] 
    table[x=scalar, y=spectral_radius, col sep=comma] {tikz_Schaubilder/src_c/_new_results/Dirichlet_MIP_DSA_1.0.csv};

    \nextgroupplot[
        legend style={
            at={($(0,0)+(1cm,1cm)$)},
            legend columns=7,
            fill=none,
            draw=black,
            anchor=center,
            align=center},
        legend to name=c_legend,
    ]
    
    \addlegendimage{empty legend}
    \addlegendentry{\textbf{Scattering Ratio:}\hspace{0.2cm}}
    
    \addplot[color={rgb,255:red,0;green,114;blue,178}, mark options={black, scale=0.75}] 
    table[x=scalar, y=spectral_radius, col sep=comma] {tikz_Schaubilder/src_c/_new_results/Robin_DSA_0.8.csv};
    \addlegendentry{0.8};

    \addplot[color={rgb,255:red,213;green,94;blue,0}, mark options={black, scale=0.75}] 
    table[x=scalar, y=spectral_radius, col sep=comma] {tikz_Schaubilder/src_c/_new_results/Robin_DSA_0.9.csv};
    \addlegendentry{0.9};
    
    \addplot[color={rgb,255:red,0;green,158;blue,115}, mark options={black, scale=0.75}] 
    table[x=scalar, y=spectral_radius, col sep=comma] {tikz_Schaubilder/src_c/_new_results/Robin_DSA_0.99.csv};
    \addlegendentry{0.99};
    
    \addplot[color={rgb,255:red,230;green,159;blue,0}, mark options={black, scale=0.75}] 
    table[x=scalar, y=spectral_radius, col sep=comma] {tikz_Schaubilder/src_c/_new_results/Robin_DSA_0.999.csv};
    \addlegendentry{0.999};
    
    \addplot[color={rgb,255:red,86;green,180;blue,233}, mark options={black, scale=0.75}] 
    table[x=scalar, y=spectral_radius, col sep=comma] {tikz_Schaubilder/src_c/_new_results/Robin_DSA_0.9999.csv};
    \addlegendentry{0.9999};

    \addplot[color={rgb,255:red,204;green,121;blue,167}, mark options={black, scale=0.75}] 
    table[x=scalar, y=spectral_radius, col sep=comma] {tikz_Schaubilder/src_c/_new_results/Robin_DSA_1.0.csv};
    \addlegendentry{1.0};

    \nextgroupplot
    
    \addlegendimage{empty legend}
    \addlegendentry{\textbf{Robin}}
    \addlegendimage{empty legend}
    \addlegendentry{\textbf{MIP}}
        
    \addplot[color={rgb,255:red,0;green,114;blue,178}, mark options={black, scale=0.75}] 
    table[x=scalar, y=spectral_radius, col sep=comma] {tikz_Schaubilder/src_c/_new_results/Robin_MIP_DSA_0.8.csv};

    \addplot[color={rgb,255:red,213;green,94;blue,0}, mark options={black, scale=0.75}] 
    table[x=scalar, y=spectral_radius, col sep=comma] {tikz_Schaubilder/src_c/_new_results/Robin_MIP_DSA_0.9.csv};
    
    \addplot[color={rgb,255:red,0;green,158;blue,115}, mark options={black, scale=0.75}] 
    table[x=scalar, y=spectral_radius, col sep=comma] {tikz_Schaubilder/src_c/_new_results/Robin_MIP_DSA_0.99.csv};
    
    \addplot[color={rgb,255:red,230;green,159;blue,0}, mark options={black, scale=0.75}] 
    table[x=scalar, y=spectral_radius, col sep=comma] {tikz_Schaubilder/src_c/_new_results/Robin_MIP_DSA_0.999.csv};
    
    \addplot[color={rgb,255:red,86;green,180;blue,233}, mark options={black, scale=0.75}] 
    table[x=scalar, y=spectral_radius, col sep=comma] {tikz_Schaubilder/src_c/_new_results/Robin_MIP_DSA_0.9999.csv};

    \addplot[color={rgb,255:red,204;green,121;blue,167}, mark options={black, scale=0.75}] 
    table[x=scalar, y=spectral_radius, col sep=comma] {tikz_Schaubilder/src_c/_new_results/Robin_MIP_DSA_1.0.csv};

    \end{groupplot}

    \node[
        anchor=north west,
        align=center,
        font=\bfseries,
        inner sep=5pt,
        fill=white,
        draw=black,
        very thick,
    ] at ([xshift=4pt,yshift=-4pt]group c1r2.north west) {Robin\\SIP};

    \node[below=1em] at (current bounding box.south)
    {\pgfplotslegendfromname{c_legend}};
    
\end{tikzpicture}
    \caption{Experiment~\ref{subsec:exp3_scattering_ratio}. Empirical
      convergence factor as a function of $\sigma_t$ for scattering
      ratios $c\in\{0.8,0.9,0.99,0.999,0.9999,1\}$. Shown are
      unaccelerated source iteration and the four DSA variants
      (SIP/MIP combined with Dirichlet/Marshak diffusion boundary
      conditions).}
    \label{fig:exp3_c_comparison}
\end{figure}

For clarity we consider spatially homogeneous scattering ratios
$c\in(0,1]$ (constant over $\Omega$), with $\sigma_s=c\sigma_t$ and
$\sigma_a=(1-c)\sigma_t$. For each fixed $c$ we vary
$\sigma_t\in[10^{-3},10^{6}]$ on a logarithmic scale, thereby probing
transport-dominated, intermediate and optically thick regimes while
isolating the effect of $c$. We use the manufactured solution
\eqref{eqn:example_solution}. Figure~\ref{fig:exp3_c_comparison}
reports the observed convergence factors for
$c\in\{0.8,0.9,0.99,0.999,0.9999,1\}$. We include $c=1$ as a stress
test. In this case the standard absorption-based contraction argument
for source iteration no longer applies, although on a bounded domain
with leakage one may still observe convergence in practice. The DSA
correction remains well defined (subject to the chosen diffusion
boundary treatment) and can still reduce the observed convergence
factor in some regimes.

Several trends are visible in Figure~\ref{fig:exp3_c_comparison}. For
the SIP variants, divergence again occurs in the intermediate regime,
while the MIP variants remain robust. Up to the onset of SIP
divergence, SIP and MIP exhibit comparable behaviour. Increasing $c$
increases the observed convergence factors, and as $c\to 1$ these
factors approach an apparent upper bound.

For the SIP variants (restricting attention to cases where they
converge), the observed upper bounds are approximately $0.56$ and $0.63$, with
average convergence factors of $0.27$ and $0.32$ for the Dirichlet and
Marshak (Robin) cases respectively. The spread is nonetheless large,
with observed values ranging from $6.1\times 10^{-4}$ up to $0.63$.

For the MIP variants (excluding the $c=1$ stress test), similar
qualitative behaviour is observed. Larger $c$ and larger $\sigma_t$
cause the observed convergence factors to remain elevated for longer
before eventually decreasing to small values (for example
$\rho\lesssim 10^{-3}$ in the most diffusive settings tested). The
maximal convergence factors observed for the Dirichlet and Robin MIP
schemes are $0.56$ and $0.62$ respectively (attained here at
$c=0.9999$), with average convergence factors of $0.15$ and $0.24$.
When $c=1$, the diffusion correction corresponds to a pure diffusion
operator
\[
\mathcal{D}=-\nabla\cdot(D\nabla\cdot), \qquad D=\frac{1}{d\sigma_t}.
\]
In this case we still observe contraction for smaller values of
$\sigma_t$ under both diffusion boundary treatments. As $\sigma_t$
increases (and therefore $D$ decreases), the diffusion correction
becomes progressively weaker and the convergence factors increase
towards their large-$\sigma_t$ limit. In the regimes where the
diffusion solve remains effective, the maximal observed convergence
factors are $0.56$ and $0.63$ for the Dirichlet and Robin MIP schemes
respectively (equal to the SIP schemes). Overall, whenever the accelerated scheme remains
convergent, MIP-based DSA provides a substantial reduction in the
observed convergence factor relative to source iteration across
high-scattering regimes.

To support these observations with representative runtime data,
Table~\ref{tab:exp3_scattering_ratios} reports wall-clock times for
selected $(c,\sigma_t)$ pairs. We include iteration counts and, where
applicable, the runtime ratio relative to source iteration. For this
table only, we increase the iteration limit to 1500 and declare
\emph{did not converge} (D.N.C.) if the iteration does not satisfy the
convergence criterion within this limit.

\begin{table}[htbp]
    \centering
    \begin{tabular}{|c|c|c c c c c|}
    \hline
    $c$ & $\sigma_t$ & Source Iteration & Dirichlet SIP & Dirichlet MIP & Robin SIP & Robin MIP  \\
    \hline
    \hline
    \multirow{3}{2em}{0.8} & 0.1 & 0.99 [18] & 2.20 [14, 44.86] & 2.13 [14, 46.35] & 2.35 [15, 42.05] & 2.27 [15, 43.53]\\
    & 1.0 & 2.97 [53] & 3.03 [21, 98.09] & 3.01 [21, 98.72] & 2.43 [17, 122.34] & 2.44 [17, 121.89]\\
    & 10 & 4.66 [86] & 3.30 [23, 141.34] & 3.27 [23, 142.74] & 2.71 [19, 172.20] & 2.71 [19, 172.28]\\
    \hline
    \multirow{3}{2em}{0.99} & 0.1 & 1.07 [18] & 2.39 [16, 44.88] & 2.40 [16, 44.83] & 2.83 [19, 37.98] & 3.00 [19, 35.84] \\
    & 1.0 & 9.54 [177] & 4.31 [30, 221.32] & 4.29 [30, 222.34] & 6.00 [42, 158.90] & 6.18 [42, 154.41] \\
    & 10 & 60.10 [1164] & 5.50 [39, 1093.37] & 5.50 [39, 1093.76] & 4.19 [30, 1433.97] & 4.21 [30, 1427.92] \\
    \hline
    \multirow{3}{2em}{0.999} & 0.1 & 1.24 [18] & 2.40 [16, 51.74] & 2.49 [16, 49.78] & 2.66 [18, 46.56] & 2.62 [18, 47.26] \\
    & 1.0 & 11.80 [224] & 4.57 [32, 257.94] & 4.69 [32, 251.71] & 7.39 [52, 159.55] & 7.37 [52, 160.17] \\
    & 10 & 76.68 [D.N.C] & 6.02 [42, -] & 5.87 [42, -] & 5.75 [41, -] & 5.72 [41, -] \\
    \hline
    \multirow{3}{2em}{1.0} & 0.1 & 1.03 [18] & 2.37 [16, 43.54] & 2.42 [16, 42.58] & 2.78 [19, 36.99] & 2.91 [19, 35.35] \\
    & 1.0 & 11.86 [227] & 4.57 [32, 259.54] & 4.64 [32, 255.36] & 7.08 [50, 167.49] & 7.11 [50, 166.81] \\
    & 10 & 76.50 [D.N.C] & 6.19 [44, -] & 6.14 [44, -] & 6.43 [46, -] & 6.58 [46, -] \\
    \hline
    \end{tabular}
    \caption{Experiment~\ref{subsec:exp3_scattering_ratio}.
      Representative wall-clock times for convergence at selected
      $(c,\sigma_t)$. Entries are reported as time [iteration count,
      runtime ratio], where the runtime ratio is
      $100\times T_{\mathrm{SI}}/T_{\mathrm{method}}$ when available.
      D.N.C.\ denotes did not converge within 1500 iterations and
      ``--'' denotes not applicable.}
    \label{tab:exp3_scattering_ratios}
\end{table}

Table~\ref{tab:exp3_scattering_ratios} confirms that iteration counts
increase as $c$ increases, most dramatically for unaccelerated source
iteration, which in the highly scattering cases shown here fails to
satisfy the stopping criterion within the iteration limit. The
accelerated schemes also require more iterations as $c$ increases, but
remain bounded in the cases shown (for example reaching $46$ iterations
for $c=1$ and $\sigma_t=10$ in this table). Overall, these results
support the conclusion that, whenever the accelerated scheme remains
convergent, DSA is effective across high-scattering regimes and
increasing optical thickness.

\subsection{Experiment 4: Changing the angular quadrature rule}
\label{subsec:exp4_quadrature}

It is well known that the convergence of source iteration and the
effectiveness of DSA depend on the angular discretisation. In
particular, the damping of error modes depends on the discrete
ordinates set, and in idealised settings (for example infinite or
periodic media) one can derive mode-by-mode damping factors that are
explicit functions of the quadrature \cite{dsa_lifeline}. On a bounded
polytopal domain, however, the convergence factor measured here also
reflects boundary effects and the interaction between the transport
sweep and the mesh geometry. 

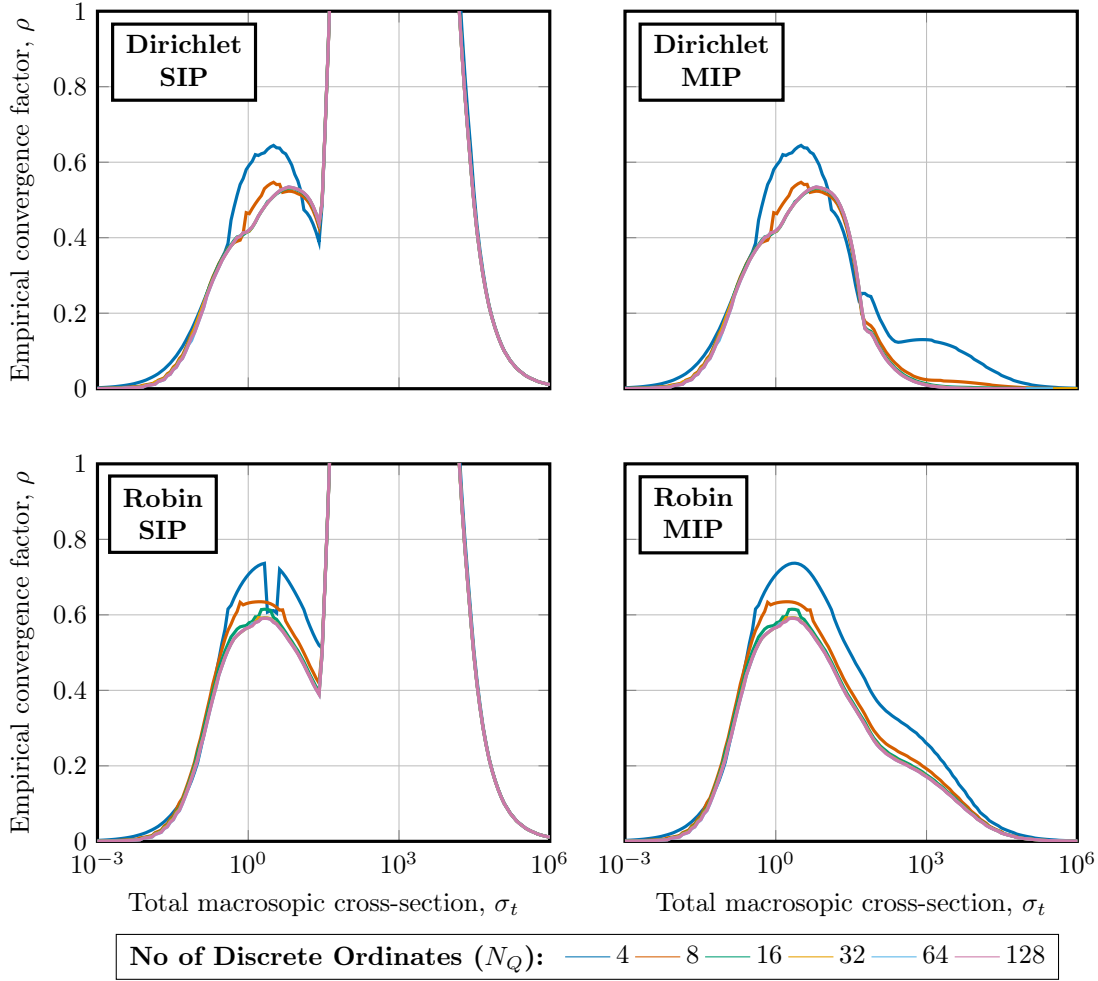
\begin{figure}[htbp]
    \centering
    \begin{tikzpicture}
    \begin{groupplot}[
        group style={
            group size=2 by 2,
            horizontal sep= 1.0cm,
            vertical sep=1cm,
            x descriptions at=edge bottom,
            y descriptions at=edge left,
            every plot/.style={
            very thick,
            ymin=0,ymax=1, 
            xmin=1e-3, xmax=1e6,
            xtick={1e-3, 1e0, 1e3, 1e6},
            legend pos=north west,
            },
        },
        enlargelimits=false,
        xlabel={Total macrosopic cross-section, $\sigma_t$},
        ylabel={Empirical convergence factor, $\rho$},
    ]
 
    % Top plot
    \nextgroupplot[
        width=0.46\textwidth,
        height=0.4\textwidth,
        xmode=log,
        grid=major,
    ]
    \addlegendimage{empty legend}
    \addlegendentry{\textbf{Dirichlet}}
    \addlegendimage{empty legend}
    \addlegendentry{\textbf{SIP}}
    
    \addplot[color={rgb,255:red,0;green,114;blue,178}, mark options={black, scale=0.75}] 
    table[x=scalar, y=spectral_radius, col sep=comma] {tikz_Schaubilder/src_q/_new_results/Dirichlet_DSA_4.csv};

    \addplot[color={rgb,255:red,213;green,94;blue,0}, mark options={black, scale=0.75}] 
    table[x=scalar, y=spectral_radius, col sep=comma] {tikz_Schaubilder/src_q/_new_results/Dirichlet_DSA_8.csv};
    
    \addplot[color={rgb,255:red,0;green,158;blue,115}, mark options={black, scale=0.75}] 
    table[x=scalar, y=spectral_radius, col sep=comma] {tikz_Schaubilder/src_q/_new_results/Dirichlet_DSA_16.csv};

    \addplot[color={rgb,255:red,230;green,159;blue,0}, mark options={black, scale=0.75}] 
    table[x=scalar, y=spectral_radius, col sep=comma] {tikz_Schaubilder/src_q/_new_results/Dirichlet_DSA_32.csv};

    \addplot[color={rgb,255:red,86;green,180;blue,233}, mark options={black, scale=0.75}] 
    table[x=scalar, y=spectral_radius, col sep=comma] {tikz_Schaubilder/src_q/_new_results/Dirichlet_DSA_64.csv};

    \addplot[color={rgb,255:red,204;green,121;blue,167}, mark options={black, scale=0.75}] 
    table[x=scalar, y=spectral_radius, col sep=comma] {tikz_Schaubilder/src_q/_new_results/Dirichlet_DSA_128.csv};
    
    \nextgroupplot[
        width=0.46\textwidth,
        height=0.4\textwidth,
        xmode=log,
        grid=major,
    ]
    \addlegendimage{empty legend}
    \addlegendentry{\textbf{Dirichlet}}
    \addlegendimage{empty legend}
    \addlegendentry{\textbf{MIP}}
    
    \addplot[color={rgb,255:red,0;green,114;blue,178}, mark options={black, scale=0.75}] 
    table[x=scalar, y=spectral_radius, col sep=comma] {tikz_Schaubilder/src_q/_new_results/Dirichlet_MIP_DSA_4.csv};

    \addplot[color={rgb,255:red,213;green,94;blue,0}, mark options={black, scale=0.75}] 
    table[x=scalar, y=spectral_radius, col sep=comma] {tikz_Schaubilder/src_q/_new_results/Dirichlet_MIP_DSA_8.csv};
    
    \addplot[color={rgb,255:red,0;green,158;blue,115}, mark options={black, scale=0.75}] 
    table[x=scalar, y=spectral_radius, col sep=comma] {tikz_Schaubilder/src_q/_new_results/Dirichlet_MIP_DSA_16.csv};

    \addplot[color={rgb,255:red,230;green,159;blue,0}, mark options={black, scale=0.75}] 
    table[x=scalar, y=spectral_radius, col sep=comma] {tikz_Schaubilder/src_q/_new_results/Dirichlet_MIP_DSA_32.csv};

    \addplot[color={rgb,255:red,86;green,180;blue,233}, mark options={black, scale=0.75}] 
    table[x=scalar, y=spectral_radius, col sep=comma] {tikz_Schaubilder/src_q/_new_results/Dirichlet_MIP_DSA_64.csv};
    
    \addplot[color={rgb,255:red,204;green,121;blue,167}, mark options={black, scale=0.75}] 
    table[x=scalar, y=spectral_radius, col sep=comma] {tikz_Schaubilder/src_q/_new_results/Dirichlet_MIP_DSA_128.csv};

    \nextgroupplot[
        width=0.46\textwidth,
        height=0.4\textwidth,
        xmode=log,
        grid=major,
        legend style={
            at={($(0,0)+(1cm,1cm)$)},
            legend columns=7,
            fill=none,
            draw=black,
            anchor=center,
            align=center},
        legend to name=legend
    ]

    % \addlegendimage{empty legend}
    % \addlegendentry{\textbf{Robin}}
    % \addlegendimage{empty legend}
    % \addlegendentry{\textbf{SIP}}
    
    \coordinate (c1) at (rel axis cs:0,1);
    \addlegendimage{empty legend}
    \addlegendentry{\textbf{No of Discrete Ordinates ($N_Q$):}\hspace{0.2cm}}
    
    \addplot[color={rgb,255:red,0;green,114;blue,178}, mark options={black, scale=0.75}] 
    table[x=scalar, y=spectral_radius, col sep=comma] {tikz_Schaubilder/src_q/_new_results/Robin_DSA_4.csv};
    \addlegendentry{4};

    \addplot[color={rgb,255:red,213;green,94;blue,0}, mark options={black, scale=0.75}] 
    table[x=scalar, y=spectral_radius, col sep=comma] {tikz_Schaubilder/src_q/_new_results/Robin_DSA_8.csv};
    \addlegendentry{8};
    
    \addplot[color={rgb,255:red,0;green,158;blue,115}, mark options={black, scale=0.75}] 
    table[x=scalar, y=spectral_radius, col sep=comma] {tikz_Schaubilder/src_q/_new_results/Robin_DSA_16.csv};
    \addlegendentry{16};

    \addplot[color={rgb,255:red,230;green,159;blue,0}, mark options={black, scale=0.75}] 
    table[x=scalar, y=spectral_radius, col sep=comma] {tikz_Schaubilder/src_q/_new_results/Robin_DSA_32.csv};
    \addlegendentry{32};

    \addplot[color={rgb,255:red,86;green,180;blue,233}, mark options={black, scale=0.75}] 
    table[x=scalar, y=spectral_radius, col sep=comma] {tikz_Schaubilder/src_q/_new_results/Robin_DSA_64.csv};
    \addlegendentry{64};
    
    \addplot[color={rgb,255:red,204;green,121;blue,167}, mark options={black, scale=0.75}] 
    table[x=scalar, y=spectral_radius, col sep=comma] {tikz_Schaubilder/src_q/_new_results/Robin_DSA_128.csv};
    \addlegendentry{128};

    \nextgroupplot[
        width=0.46\textwidth,
        height=0.4\textwidth,
        xmode=log,
        grid=major,
    ]
    \addlegendimage{empty legend}
    \addlegendentry{\textbf{Robin}}
    \addlegendimage{empty legend}
    \addlegendentry{\textbf{MIP}}
    
    \coordinate (c2) at (rel axis cs:1,1);
    
    \addplot[color={rgb,255:red,0;green,114;blue,178}, mark options={black, scale=0.75}] 
    table[x=scalar, y=spectral_radius, col sep=comma] {tikz_Schaubilder/src_q/_new_results/Robin_MIP_DSA_4.csv};

    \addplot[color={rgb,255:red,213;green,94;blue,0}, mark options={black, scale=0.75}] 
    table[x=scalar, y=spectral_radius, col sep=comma] {tikz_Schaubilder/src_q/_new_results/Robin_MIP_DSA_8.csv};
    
    \addplot[color={rgb,255:red,0;green,158;blue,115}, mark options={black, scale=0.75}] 
    table[x=scalar, y=spectral_radius, col sep=comma] {tikz_Schaubilder/src_q/_new_results/Robin_MIP_DSA_16.csv};

    \addplot[color={rgb,255:red,230;green,159;blue,0}, mark options={black, scale=0.75}] 
    table[x=scalar, y=spectral_radius, col sep=comma] {tikz_Schaubilder/src_q/_new_results/Robin_MIP_DSA_32.csv};

    \addplot[color={rgb,255:red,86;green,180;blue,233}, mark options={black, scale=0.75}] 
    table[x=scalar, y=spectral_radius, col sep=comma] {tikz_Schaubilder/src_q/_new_results/Robin_MIP_DSA_64.csv};
    
    \addplot[color={rgb,255:red,204;green,121;blue,167}, mark options={black, scale=0.75}] 
    table[x=scalar, y=spectral_radius, col sep=comma] {tikz_Schaubilder/src_q/_new_results/Robin_MIP_DSA_128.csv};

    \end{groupplot}

    \node[
        anchor=north west,
        align=center,
        font=\bfseries,
        inner sep=5pt,
        fill=white,
        draw=black,
        very thick,
    ] at ([xshift=4pt,yshift=-4pt]group c1r2.north west) {Robin\\SIP};

    \coordinate (c3) at ($(c1)!.5!(c2)$);
    \node[below] at (c3 |- current bounding box.south)
      {\pgfplotslegendfromname{legend}};
    
\end{tikzpicture}
    \caption{Experiment~\ref{subsec:exp4_quadrature}. Empirical
      convergence factor as a function of $\sigma_t$ for angular
      refinements $N_Q\in\{4,8,16,32,64,128\}$, comparing the four DSA
      variants (SIP/MIP combined with Dirichlet/Marshak diffusion
      boundary conditions).}
    \label{fig:exp4_q_comparison}
\end{figure}

In this experiment we refine the angular discretisation by increasing
the number of discrete ordinates. We consider
\[
N_Q\in\{4,8,16,32,64,128\},
\]
and, for each $N_Q$, compute the observed convergence factor $\rho$
over the same range of $\sigma_t$ values as in
Experiment~\ref{subsec:exp1_bc_penalty}. The results are shown in
Figure~\ref{fig:exp4_q_comparison}.

Figure~\ref{fig:exp4_q_comparison} shows that, for each scheme, the
observed convergence factors approach a limiting curve as $N_Q$
increases. For the SIP variants divergence occurs over essentially the
same range of $\sigma_t$ values, largely independent of $N_Q$. In
particular, refining the angular discretisation does not materially
shift the range of $\sigma_t$ over which the SIP-based DSA fails.

For $N_Q\gtrsim 32$, further angular refinement yields only small
changes in the measured convergence factors, the curves are visually
close for $N_Q=32,64,128$. For the MIP variants, the convergence
remains well behaved across the full range of $\sigma_t$ for every
$N_Q$, again approaching a limiting curve as $N_Q$ increases. In the
MIP--Dirichlet case we observe a noticeable steep drop in $\rho$ in the
regime $\sigma_t\sim 10^{2}$, whereas for MIP--Marshak the behaviour
in this regime is less pronounced and is instead reflected by a
slower decrease in $\rho$ as $\sigma_t$ increases. As in
Experiments~\ref{subsec:exp1_bc_penalty} and
\ref{subsec:exp3_scattering_ratio}, these differences occur in the
transition between transport-dominated and diffusion-dominated
behaviour, where the boundary model and the jump control most strongly
influence the correction.

The iteration counts required for convergence are broadly consistent
across $N_Q$, particularly for $N_Q\ge 16$. The largest iteration
counts in this experiment occur when $\rho\approx 0.6$ for the MIP--Marshak
variant, reaching approximately 55 iterations.

Finally, we remark that very small ordinate sets can exhibit classical
ray effects, which can distort angular error modes and thereby
influence the measured convergence rates. While such artefacts were
not visually dominant here, the smallest quadrature sets do not
resolve the angular dependence of \eqref{eqn:example_solution} well,
which is consistent with the markedly different convergence factors
observed for $N_Q\in\{4,8\}$ in
Figure~\ref{fig:exp4_q_comparison}.

\subsection{Experiment 5: Effect of refining the mesh}
\label{subsec:exp5_h_refinement}

We now consider how spatial refinement ($h$-refinement) affects the
effectiveness of DSA. Refinement is a standard practice in DG methods
and it is therefore important to understand how the convergence
properties of the accelerated iteration behave as the mesh is refined.

\begin{figure}[htbp]
  \centering
  \includegraphics[width=\textwidth]{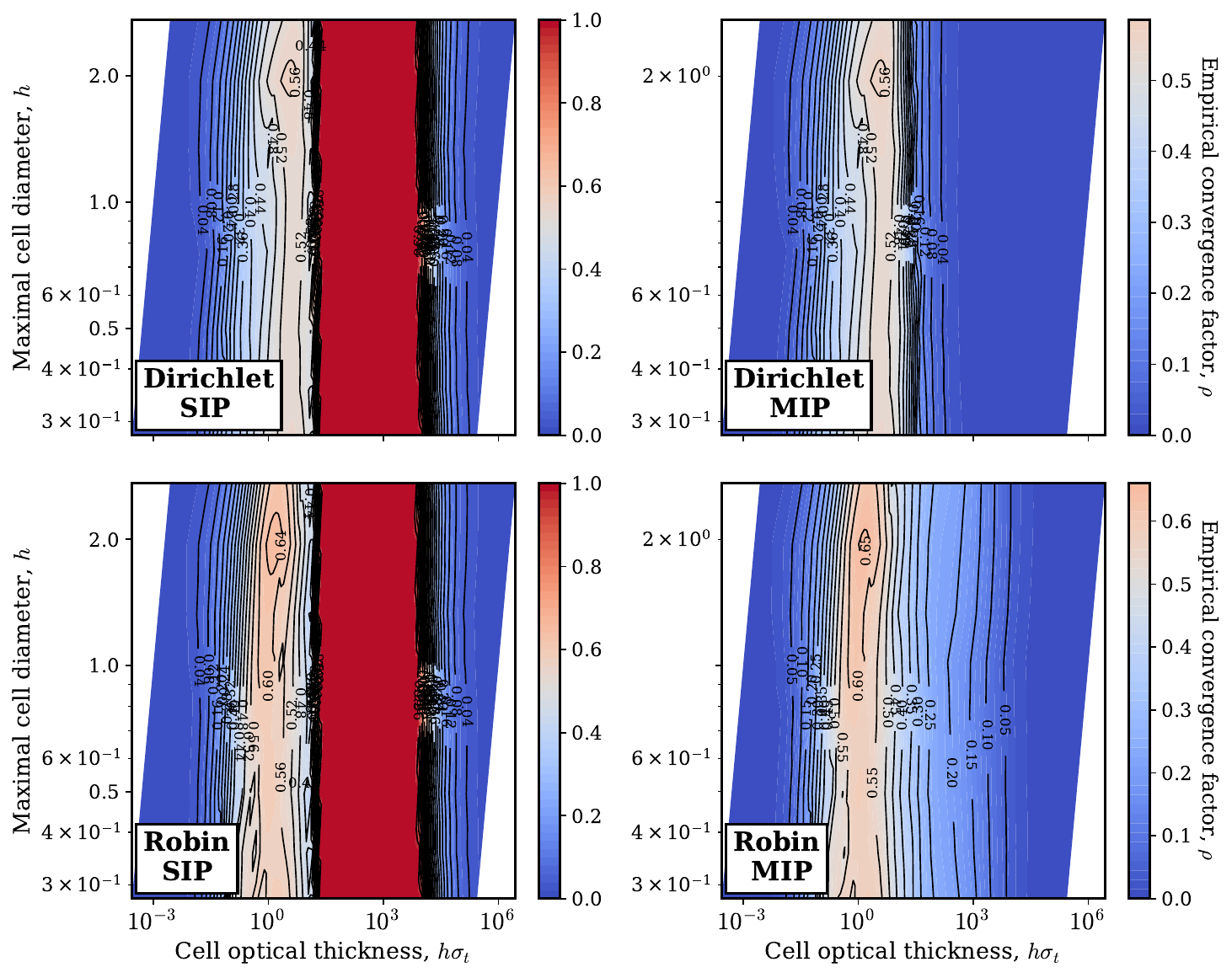}
  \caption{Experiment~\ref{subsec:exp5_h_refinement}. Empirical
    convergence factor as a function of the optical thickness
    $h\sigma_t$, comparing unaccelerated source iteration with the
    four DSA variants (SIP/MIP combined with Dirichlet/Marshak
    diffusion boundary conditions) across meshes of increasing
    refinement.}
  \label{fig:exp5_h_comparison}
\end{figure}

We refine globally by re-meshing the domain with increasing numbers of
elements. In practice, refinement is often localised, but since
conditioning and penalty scaling are primarily influenced by the
smallest length scales, a global study provides a clean
baseline. Specifically, we remesh $\Omega=(0,10)^2$ for the
manufactured solution \eqref{eqn:example_solution} using
\[
|\mathcal{T}_h|\in\{32,64,128,256,512,1024,2048,4096\},
\]
and take the corresponding values of
$h=\max_{\kappa\in\mathcal{T}_h}h_\kappa$. To isolate mesh-size
effects from anisotropy effects, each mesh is smoothed using 10 Lloyd
iterations. Since the DSA equation is solved on the same mesh and the
interior penalty depends on both geometry and $\sigma_t$, the
diffusion operator $\vec{D}$ is reassembled for each refinement
level. Results are reported in terms of the optical thickness
$h\sigma_t$, which is the natural dimensionless parameter governing
the transition between transport and diffusion
regimes. Figure~\ref{fig:exp5_h_comparison} shows the resulting
empirical convergence factors as functions of $h\sigma_t$.

Figure~\ref{fig:exp5_h_comparison} indicates that, across
refinements, the observed convergence factor largely collapses when
plotted against the optical thickness $h\sigma_t$, and this appears to
hold for both SIP- and MIP-based corrections. In the MIP variants,
decreasing $h$ leads to a modest but consistent reduction in the
observed convergence factor in the intermediate regime. This is
encouraging for highly refined or complex meshes, including
large-scale simulations with very large element counts.

Across the refinement levels tested, typical maximal convergence
factors are approximately $0.57$ and $0.66$ for the Dirichlet and
Marshak SIP schemes respectively (in regimes where the SIP variants
remain convergent), and the same for both MIP schemes. For
the MIP variants these peak values correspond to around $48$ and $66$ iterations
to convergence respectively, with a slight decrease in iteration count as the mesh is refined.

\begin{figure}[htbp]
    \centering
    \begin{tikzpicture}
    \begin{groupplot}[
        group style={
            group size=2 by 2,
            horizontal sep= 1.0cm,
            vertical sep=1cm,
            x descriptions at=edge bottom,
            y descriptions at=edge left,
            every plot/.style={
            very thick,
            ymin=0,ymax=70, 
            xtick={1e-3, 1e0, 1e3, 1e6},
            legend pos=north west,
            },
        },
        enlargelimits=false,
        xlabel={Cell optical thickness, $h\sigma_t$},
        ylabel={Number of iterations},
    ]
 
    % Top plot
    \nextgroupplot[
        width=0.46\textwidth,
        height=0.4\textwidth,
        xmode=log,
        grid=major,
    ]
    \addlegendimage{empty legend}
    \addlegendentry{\textbf{Dirichlet}}
    \addlegendimage{empty legend}
    \addlegendentry{\textbf{SIP}}
    
    \addplot[color={rgb,255:red,0;green,114;blue,178}, mark options={black, scale=0.75}] 
    table[x expr= \thisrow{scalar} * \thisrow{h}, y=n_iterations, col sep=comma] {tikz_Schaubilder/src_h/_new_results/Dirichlet_DSA_64.csv};

    \addplot[color={rgb,255:red,213;green,94;blue,0}, mark options={black, scale=0.75}] 
    table[x expr= \thisrow{scalar} * \thisrow{h}, y=n_iterations, col sep=comma] {tikz_Schaubilder/src_h/_new_results/Dirichlet_DSA_512.csv};
    
    \addplot[color={rgb,255:red,0;green,158;blue,115}, mark options={black, scale=0.75}] 
    table[x expr= \thisrow{scalar} * \thisrow{h}, y=n_iterations, col sep=comma] {tikz_Schaubilder/src_h/_new_results/Dirichlet_DSA_4096.csv};
    
    \nextgroupplot[
        width=0.46\textwidth,
        height=0.4\textwidth,
        xmode=log,
        grid=major,
    ]
    \addlegendimage{empty legend}
    \addlegendentry{\textbf{Dirichlet}}
    \addlegendimage{empty legend}
    \addlegendentry{\textbf{MIP}}
    
    \addplot[color={rgb,255:red,0;green,114;blue,178}, mark options={black, scale=0.75}] 
    table[x expr= \thisrow{scalar} * \thisrow{h}, y=n_iterations, col sep=comma] {tikz_Schaubilder/src_h/_new_results/Dirichlet_MIP_DSA_64.csv};

    \addplot[color={rgb,255:red,213;green,94;blue,0}, mark options={black, scale=0.75}] 
    table[x expr= \thisrow{scalar} * \thisrow{h}, y=n_iterations, col sep=comma] {tikz_Schaubilder/src_h/_new_results/Dirichlet_MIP_DSA_512.csv};
    
    \addplot[color={rgb,255:red,0;green,158;blue,115}, mark options={black, scale=0.75}] 
    table[x expr= \thisrow{scalar} * \thisrow{h}, y=n_iterations, col sep=comma] {tikz_Schaubilder/src_h/_new_results/Dirichlet_MIP_DSA_4096.csv};

    \nextgroupplot[
        width=0.46\textwidth,
        height=0.4\textwidth,
        xmode=log,
        grid=major,
        legend style={
            at={($(0,0)+(1cm,1cm)$)},
            legend columns=7,
            fill=none,
            draw=black,
            anchor=center,
            align=center},
        legend to name=i_legend
    ]
    
    \coordinate (c1) at (rel axis cs:0,1);
    \addlegendimage{empty legend}
    \addlegendentry{\textbf{No. Elements:}\hspace{0.2cm}}
    
    \addplot[color={rgb,255:red,0;green,114;blue,178}, mark options={black, scale=0.75}] 
    table[x expr= \thisrow{scalar} * \thisrow{h}, y=n_iterations, col sep=comma] {tikz_Schaubilder/src_h/_new_results/Robin_DSA_64.csv};
    \addlegendentry{64};

    \addplot[color={rgb,255:red,213;green,94;blue,0}, mark options={black, scale=0.75}] 
    table[x expr= \thisrow{scalar} * \thisrow{h}, y=n_iterations, col sep=comma] {tikz_Schaubilder/src_h/_new_results/Robin_DSA_512.csv};
    \addlegendentry{512};
    
    \addplot[color={rgb,255:red,0;green,158;blue,115}, mark options={black, scale=0.75}] 
    table[x expr= \thisrow{scalar} * \thisrow{h}, y=n_iterations, col sep=comma] {tikz_Schaubilder/src_h/_new_results/Robin_DSA_4096.csv};
    \addlegendentry{4096};

    \nextgroupplot[
        width=0.46\textwidth,
        height=0.4\textwidth,
        xmode=log,
        grid=major,
    ]
    \addlegendimage{empty legend}
    \addlegendentry{\textbf{Robin}}
    \addlegendimage{empty legend}
    \addlegendentry{\textbf{MIP}}
    
    \coordinate (c2) at (rel axis cs:1,1);
    
    \addplot[color={rgb,255:red,0;green,114;blue,178}, mark options={black, scale=0.75}] 
    table[x expr= \thisrow{scalar} * \thisrow{h}, y=n_iterations, col sep=comma] {tikz_Schaubilder/src_h/_new_results/Robin_MIP_DSA_64.csv};

    \addplot[color={rgb,255:red,213;green,94;blue,0}, mark options={black, scale=0.75}] 
    table[x expr= \thisrow{scalar} * \thisrow{h}, y=n_iterations, col sep=comma] {tikz_Schaubilder/src_h/_new_results/Robin_MIP_DSA_512.csv};
    
    \addplot[color={rgb,255:red,0;green,158;blue,115}, mark options={black, scale=0.75}] 
    table[x expr= \thisrow{scalar} * \thisrow{h}, y=n_iterations, col sep=comma] {tikz_Schaubilder/src_h/_new_results/Robin_MIP_DSA_4096.csv};

    \end{groupplot}

    \node[
        anchor=north west,
        align=center,
        font=\bfseries,
        inner sep=5pt,
        fill=white,
        draw=black,
        very thick,
    ] at ([xshift=4pt,yshift=-4pt]group c1r2.north west) {Robin\\SIP};

    \coordinate (c3) at ($(c1)!.5!(c2)$);
    \node[below] at (c3 |- current bounding box.south)
      {\pgfplotslegendfromname{i_legend}};
    
\end{tikzpicture}
    \caption{Experiment~\ref{subsec:exp5_h_refinement}. Iteration
      counts to convergence as a function of $\sigma_t$ for the four
      DSA variants, comparing multiple mesh refinement levels.}
    \label{fig:exp5_h_iterations}
\end{figure}
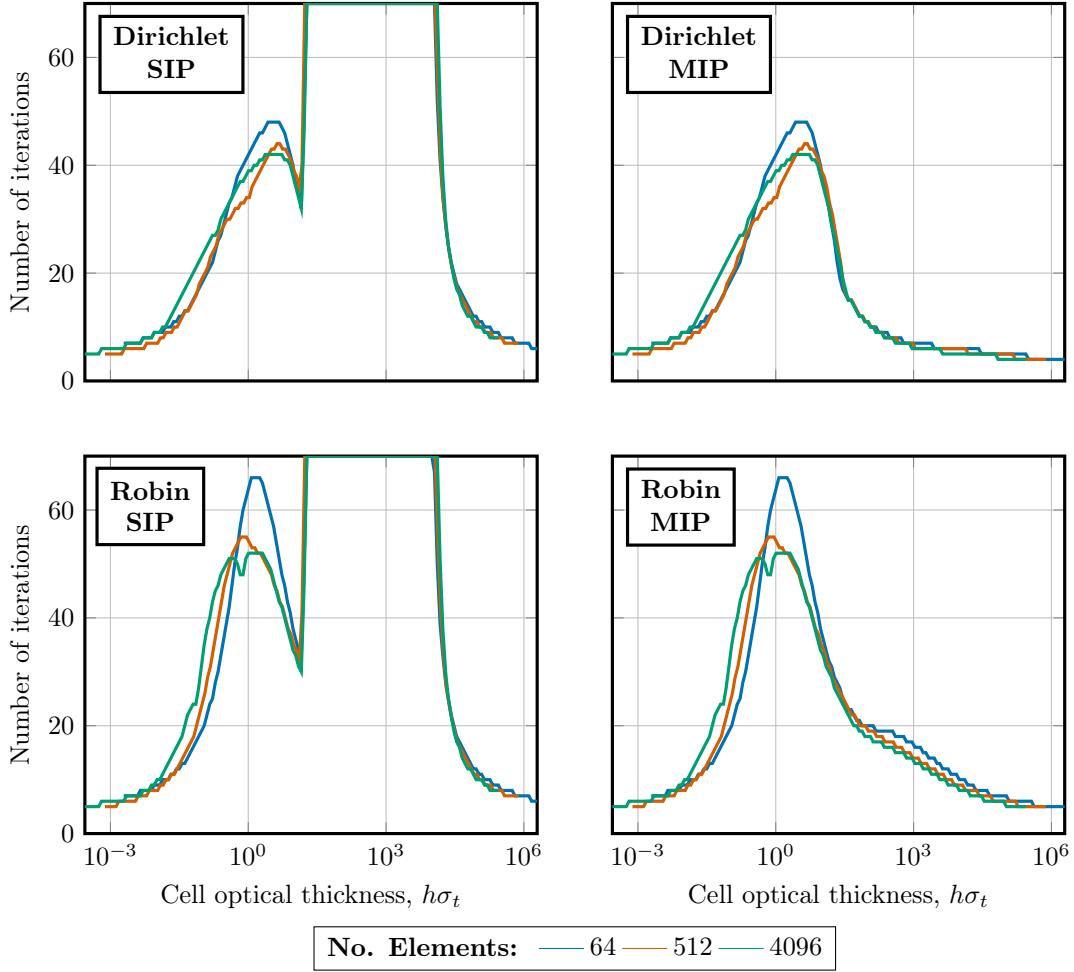

To complement Figure~\ref{fig:exp5_h_comparison},
Figure~\ref{fig:exp5_h_iterations} reports the iteration counts
required for convergence for each DSA variant across the refinement
levels.

The iteration-count trends in Figure~\ref{fig:exp5_h_iterations}
mirror the behaviour of the observed convergence factors in
Figure~\ref{fig:exp5_h_comparison}, as expected. In contrast, the
unaccelerated source iteration results in
Experiment~\ref{subsec:exp1_bc_penalty} show that iteration counts
grow rapidly with increasing optical thickness and exceed the
iteration cap early in the diffusive regime.

Differences between the Marshak and Dirichlet diffusion boundary
conditions remain pronounced under mesh refinement. The most visible
distinction occurs for the MIP variants in the intermediate
regime. The Dirichlet correction exhibits a steep drop-off in $\rho$ as
a function of $h\sigma_t$, while the Marshak correction exhibits a
longer, slower decay. A plausible explanation is the effectiveness of the boundary conditions in the intermediate regime. Potentially small inconsistencies between the asymptotically correct boundary conditions and the applied Marshak (Robin) boundary conditions may result in a lack of damping of the low-frequency boundary error components for a subset of optical thicknesses  \cite{malvagi1991initial}.
This interpretation is consistent with the scaling of the
boundary layer thickness as $\mathcal{O}(1/\sigma_t)$ and with the
fact that the two boundary treatments become increasingly similar in
the strongly diffusive regime (Figure~\ref{fig:exp5_h_comparison}).
Nonetheless, both boundary treatments approach similar upper bounds
for the observed convergence factor in the thick limit.

\subsection{Experiment 6: Effect of varying the polynomial degree}
\label{subsec:exp6_p_refinement}

While $h$-refinement is often effective for resolving geometric and
transport features, $p$-refinement is frequently the most direct route
to reducing discretisation error. It is, however, more costly the
number of local degrees of freedom per element grows like
$\binom{p+d}{d}$ and, correspondingly, the global matrices increase in
both size and fill. The interior penalty parameter
$\sigma_{\mathrm{SIP}}$ also depends on $p$ (notably through a $p^2$
scaling), so changes in $p$ are expected to affect the stability and
performance of DSA. We therefore investigate the effect of
$p$-refinement on the accelerated iteration.

In this experiment, the transport and diffusion correction problems
are discretised using the same polynomial degree space (we comment on
mixed-degree choices later). We consider $p\in\{1,2,3\}$ and use the
manufactured solution \eqref{eqn:example_solution}. The resulting
observed convergence factors as functions of $\sigma_t$ are shown in
Figure~\ref{fig:exp6_p_comparison}.

\begin{figure}[htbp]
    \centering
    \begin{tikzpicture}
    \begin{groupplot}[
        group style={
            group size=2 by 2,
            horizontal sep= 1.0cm,
            vertical sep=1cm,
            x descriptions at=edge bottom,
            y descriptions at=edge left,
            every plot/.style={
            very thick,
            ymin=0,ymax=1, 
            xmin=1e-3, xmax=1e6,
            xtick={1e-3, 1e0, 1e3, 1e6},
            legend pos=north west,
            },
        },
        enlargelimits=false,
        xlabel={Total macrosopic cross-section, $\sigma_t$},
        ylabel={Empirical convergence factor, $\rho$},
    ]
 
    % Top plot
    \nextgroupplot[
        width=0.46\textwidth,
        height=0.4\textwidth,
        xmode=log,
        grid=major,
    ]
    \addlegendimage{empty legend}
    \addlegendentry{\textbf{Dirichlet}}
    \addlegendimage{empty legend}
    \addlegendentry{\textbf{SIP}}
    
    \addplot[color={rgb,255:red,0;green,114;blue,178}, mark options={black, scale=0.75}] 
    table[x=scalar, y=spectral_radius, col sep=comma] {tikz_Schaubilder/src_p/_new_results/Dirichlet_DSA_1.csv};

    \addplot[color={rgb,255:red,213;green,94;blue,0}, mark options={black, scale=0.75}] 
    table[x=scalar, y=spectral_radius, col sep=comma] {tikz_Schaubilder/src_p/_new_results/Dirichlet_DSA_2.csv};
    
    \addplot[color={rgb,255:red,0;green,158;blue,115}, mark options={black, scale=0.75}] 
    table[x=scalar, y=spectral_radius, col sep=comma] {tikz_Schaubilder/src_p/_new_results/Dirichlet_DSA_3.csv};
    
    \nextgroupplot[
        width=0.46\textwidth,
        height=0.4\textwidth,
        xmode=log,
        grid=major,
    ]
    \addlegendimage{empty legend}
    \addlegendentry{\textbf{Dirichlet}}
    \addlegendimage{empty legend}
    \addlegendentry{\textbf{MIP}}
    
    \addplot[color={rgb,255:red,0;green,114;blue,178}, mark options={black, scale=0.75}] 
    table[x=scalar, y=spectral_radius, col sep=comma] {tikz_Schaubilder/src_p/_new_results/Dirichlet_MIP_DSA_1.csv};

    \addplot[color={rgb,255:red,213;green,94;blue,0}, mark options={black, scale=0.75}] 
    table[x=scalar, y=spectral_radius, col sep=comma] {tikz_Schaubilder/src_p/_new_results/Dirichlet_MIP_DSA_2.csv};
    
    \addplot[color={rgb,255:red,0;green,158;blue,115}, mark options={black, scale=0.75}] 
    table[x=scalar, y=spectral_radius, col sep=comma] {tikz_Schaubilder/src_p/_new_results/Dirichlet_MIP_DSA_3.csv};

    \nextgroupplot[
        width=0.46\textwidth,
        height=0.4\textwidth,
        xmode=log,
        grid=major,
        legend style={
            at={($(0,0)+(1cm,1cm)$)},
            legend columns=7,
            fill=none,
            draw=black,
            anchor=center,
            align=center},
        legend to name=p_legend
    ]

    % \addlegendimage{empty legend}
    % \addlegendentry{\textbf{Robin}}
    % \addlegendimage{empty legend}
    % \addlegendentry{\textbf{SIP}}
    
    \coordinate (c1) at (rel axis cs:0,1);
    \addlegendimage{empty legend}
    \addlegendentry{\textbf{Polynomial degree ($p$):}\hspace{0.2cm}}
    
    \addplot[color={rgb,255:red,0;green,114;blue,178}, mark options={black, scale=0.75}] 
    table[x=scalar, y=spectral_radius, col sep=comma] {tikz_Schaubilder/src_p/_new_results/Robin_DSA_1.csv};
    \addlegendentry{1};

    \addplot[color={rgb,255:red,213;green,94;blue,0}, mark options={black, scale=0.75}] 
    table[x=scalar, y=spectral_radius, col sep=comma] {tikz_Schaubilder/src_p/_new_results/Robin_DSA_2.csv};
    \addlegendentry{2};
    
    \addplot[color={rgb,255:red,0;green,158;blue,115}, mark options={black, scale=0.75}] 
    table[x=scalar, y=spectral_radius, col sep=comma] {tikz_Schaubilder/src_p/_new_results/Robin_DSA_3.csv};
    \addlegendentry{3};

    \nextgroupplot[
        width=0.46\textwidth,
        height=0.4\textwidth,
        xmode=log,
        grid=major,
    ]
    \addlegendimage{empty legend}
    \addlegendentry{\textbf{Robin}}
    \addlegendimage{empty legend}
    \addlegendentry{\textbf{MIP}}
    
    \coordinate (c2) at (rel axis cs:1,1);
    
    \addplot[color={rgb,255:red,0;green,114;blue,178}, mark options={black, scale=0.75}] 
    table[x=scalar, y=spectral_radius, col sep=comma] {tikz_Schaubilder/src_p/_new_results/Robin_MIP_DSA_1.csv};

    \addplot[color={rgb,255:red,213;green,94;blue,0}, mark options={black, scale=0.75}] 
    table[x=scalar, y=spectral_radius, col sep=comma] {tikz_Schaubilder/src_p/_new_results/Robin_MIP_DSA_2.csv};
    
    \addplot[color={rgb,255:red,0;green,158;blue,115}, mark options={black, scale=0.75}] 
    table[x=scalar, y=spectral_radius, col sep=comma] {tikz_Schaubilder/src_p/_new_results/Robin_MIP_DSA_3.csv};

    \end{groupplot}

    \node[
        anchor=north west,
        align=center,
        font=\bfseries,
        inner sep=5pt,
        fill=white,
        draw=black,
        very thick,
    ] at ([xshift=4pt,yshift=-4pt]group c1r2.north west) {Robin\\SIP};

    \coordinate (c3) at ($(c1)!.5!(c2)$);
    \node[below] at (c3 |- current bounding box.south)
      {\pgfplotslegendfromname{p_legend}};
    
\end{tikzpicture}
    \caption{Experiment~\ref{fig:exp6_p_comparison}. Empirical
      convergence factor as a function of $\sigma_t$ for polynomial
      degrees $p\in\{1,2,3\}$, comparing the four DSA variants
      (SIP/MIP combined with Dirichlet/Marshak diffusion boundary
      conditions).}
    \label{fig:exp6_p_comparison}
\end{figure}

Figure~\ref{fig:exp6_p_comparison} shows that, without MIP, divergence
is observed for $p=1,2,3$. Consistent with the dependence of
$\sigma_{\mathrm{SIP}}$ on $p$, the onset of divergence occurs at
larger values of $\sigma_t$ as $p$ increases. Interestingly,
convergent behaviour re-emerges for all polynomial degrees once
$\sigma_t$ exceeds approximately $10^4$. In this regime the reaction
term becomes dominant since $\sigma_a=(1-c)\sigma_t$ increases
proportionally with $\sigma_t$, coercivity in the diffusion problem
strengthens and the reaction term overwhelms the under-penalised
diffusion contribution.

Comparing with the MIP variants, the observed convergence factors
remain similar to the SIP variants up to the divergence point, after
which the MIP formulation maintains stable behaviour. For both
Dirichlet and Marshak diffusion boundary conditions, the MIP
convergence factors plateau and subsequently decrease towards zero as
$\sigma_t$ increases, indicating effective damping of the remaining
error modes. Crucially, the MIP variants converge for all polynomial
degrees tested across the full range of $\sigma_t$. Iteration counts
are typically less than 60 for most of the intermediate
regime $\sigma_t\in(10^{0},10^{3})$; for $p=1,2,3$, the median
convergence factors are $0.38$, $0.53$, $0.60$ for the Dirichlet MIP
variant and $0.37$, $0.58$, $0.62$ for the Marshak MIP variant,
which shows only modest variation with $p$ in challenging regimes.

The Dirichlet MIP curves also appear to converge towards one another
as $p$ increases, suggesting the existence of a limiting convergence
profile in $p$. This limiting behaviour remains well below unity,
indicating robust convergence even as conditioning deteriorates with
increasing polynomial degree. Similar trends are observed for the
Marshak MIP variant, although the behaviour is less
pronounced. Overall, these results support the use of $p$-refinement
within the present polytopal DGFEM--DSA framework.

To complement Figure~\ref{fig:exp6_p_comparison}, we examine how the
total runtime advantage provided by DSA scales with $p$ in a
representative case where unaccelerated source iteration is not
prohibitively slow. We set $\sigma_t=0.5$ and keep the remaining
parameters as in the numerical setup of
Section~\ref{section:numerics}. We remove the iteration cap for this
timing study so that all schemes run to the convergence criterion. The
results are reported in Table~\ref{tab:exp6_p_speedup}.

\begin{table}[htbp]
    \centering
    \begin{subtable}{\linewidth}
        \centering
        \begin{tabular}{||c| c c c c c||} 
         \hline
            $p$ & Source Iteration & Dirichlet SIP & Dirichlet MIP & Robin SIP & Robin MIP \\ [0.5ex] 
             \hline\hline
             1 & \hspace{0.175cm}4.92 [84] & \hspace{0.35cm}4.55   [29] & \hspace{0.35cm}4.56    [29] & \hspace{0.175cm}7.39  [47] & \hspace{0.175cm}6.92  [47] \\ 
             2 & \hspace{0.175cm}9.27 [47] & \hspace{0.35cm}8.99   [25] & \hspace{0.35cm}8.95    [25] & 13.83 [39]                 & 13.77 [39] \\
             3 & 13.46 [43]                & \hspace{0.175cm}16.99 [24] & \hspace{0.175cm}16.78  [24] & 23.56 [33]                 & 23.66 [33] \\
             4 & 17.14 [38]                & \hspace{0.175cm}42.58 [24] & \hspace{0.175cm}16.77  [24] & 47.12 [26]                 & 47.60 [26] \\
             5 & 33.10 [38]                & 116.41 [24]                & 116.61 [24]                 & 94.59 [19]                 & 93.36 [19] \\[1ex] 
             \hline
        \end{tabular}
        \vspace{0.1cm}
        \caption{Wall-clock time (s) [iterations].}
        \vspace{0.1cm}
    \end{subtable}\par
    \begin{subtable}{\linewidth}
        \centering
        \begin{tabular}{||c| c c c c||} 
         \hline
         \multicolumn{1}{||p{0.5cm}|}{\centering $p$}  & \multicolumn{1}{|p{1.5cm}|}{\centering Dirichlet \\ SIP} &  \multicolumn{1}{|p{1.5cm}|}{\centering Dirichlet \\ MIP} &  \multicolumn{1}{|p{1.0cm}|}{\centering Robin \\ SIP} & \multicolumn{1}{|p{1.0cm}||}{\centering Robin \\ MIP} \\ [0.5ex] 
         \hline\hline
         1 & 108.12 & 107.70 & \hspace{0.175cm}66.50 & \hspace{0.175cm}71.01 \\ 
         2 & 103.16 & 103.61 & \hspace{0.175cm}67.06 & \hspace{0.175cm}67.32 \\
         3 & \hspace{0.175cm}79.23  & \hspace{0.175cm}80.22  & \hspace{0.175cm}57.12 & \hspace{0.175cm}56.87 \\
         4 & \hspace{0.175cm}40.25  & \hspace{0.175cm}40.07  & \hspace{0.175cm}36.38  & \hspace{0.175cm}36.02  \\
         5 & \hspace{0.175cm}28.43  & \hspace{0.175cm}28.38  & \hspace{0.175cm}34.99  & \hspace{0.175cm}35.45  \\[1ex] 
         \hline
        \end{tabular}
        \vspace{0.1cm}
        \caption{Runtime ratio (\%) relative to source iteration, defined by $100\times T_{\mathrm{SI}}/T_{\mathrm{method}}$.}
    \end{subtable}
    \caption{Experiment~\ref{subsec:exp6_p_refinement}. Wall-clock times and runtime ratios for $\sigma_t=0.5$ as a function of the polynomial degree $p$.}
    \label{tab:exp6_p_speedup}
\end{table} 

Table~\ref{tab:exp6_p_speedup} leads to a complementary conclusion to
Figure~\ref{fig:exp6_p_comparison}. For this comparatively mild choice
$\sigma_t=0.5$, unaccelerated source iteration already has a modest
observed convergence factor, and as $p$ increases the runtime
advantage provided by DSA becomes less significant. Indeed, source
iteration can become the fastest method despite requiring more outer
iterations, because each DSA iteration incurs the additional diffusion
correction solve. This behaviour is expected whenever the convergence
of source iteration is not severely degraded. In contrast, once
$\sigma_t$ enters the intermediate and diffusive regimes (for example
$\sigma_t\gtrsim 10$), the iteration counts for source iteration grow
rapidly and acceleration becomes necessary for all polynomial degrees,
consistent with the trends in Figure~\ref{fig:exp6_p_comparison}.

\begin{remark}[Effect of using different polynomial degrees in transport and diffusion]
The DSA correction is designed to approximate the \emph{low-order,
slowly damped} components of the transport iteration error. In the
discrete setting, the correction is computed by solving a diffusion
problem whose right-hand side is formed from the scalar-flux
difference $\phi_h^{(n+\frac12)}-\phi_h^{(n)}$. If the diffusion
problem is discretised in a different space $\mathbb{V}_{p_D}$ from
the transport space $\mathbb{V}_{p_T}$, then a transfer operator is
required to map the transport difference into the diffusion space and
to map the correction back. In practice this is typically an $L^2$
projection (or a mass-lumped variant) on each element, but any such
transfer introduces an additional approximation error and can suppress
precisely the low-frequency content that DSA is intended to remove. If
$p_D<p_T$ and the diffusion space is too coarse, the correction can
become under-resolved and the accelerated iteration may revert towards
the unaccelerated convergence behaviour, or even stagnate, despite
each iteration being more expensive. Conversely, taking
$p_D>p_T$ generally increases the cost of the diffusion solve without
improving the convergence factor, since the dominant error modes being
targeted are already low order. For this reason, and to maintain the
consistency properties established in \cite{dsaanal}, we take
$p_D=p_T$ throughout the main experiments and consider mixed-degree
choices only in a dedicated study.
\end{remark}

\subsection{Experiment 7: Effect of anisotropic elements}
\label{subsec:exp7_anisotropy}

Anisotropic elements are known to present challenges for DG
discretisations and for DSA-based acceleration. In particular,
anisotropy can lead to under-penalisation and deterioration of
stability unless the interior penalty is chosen carefully. These
effects have been discussed in the transport--DSA context by Wang and
Ragusa \cite{Wang01102010}, Turcksin and Ragusa
\cite{TURCKSIN2014356}, and Ragusa \cite{RAGUSA2015195}, and in the
broader DG literature, for example by Georgoulis
\cite{georgoulis2006hp}.

Since Voronoi tessellations are attractive in transport applications
(for example in moving-mesh contexts \cite{Munoz2012} and for modern
basis constructions \cite{boscheri2024local, boscheri2022continuous}),
we assess DSA performance on meshes with potentially high anisotropy.
For a polytopal element $\kappa$, let $h_{\kappa,\max}$ and
$h_{\kappa,\min}$ denote the largest and smallest edge lengths of a
bounding box of $\kappa$ (not necessarily aligned with the Cartesian
axes). We define the global anisotropy ratio by
\begin{equation}
    \label{eqn:anisotropy_ratio}
    \eta := \max_{\kappa\in\mathcal{T}_h}\left\{\frac{h_{\kappa,\max}}{h_{\kappa,\min}}\right\}.
\end{equation}
This yields a single mesh-wise anisotropy measure. For the
Lloyd-smoothed meshes used in
Experiments~\ref{subsec:exp1_bc_penalty}--\ref{subsec:exp6_p_refinement},
we typically have $\eta<1.5$.

\begin{figure}[htbp]
  \centering
  \includegraphics[width=\textwidth]{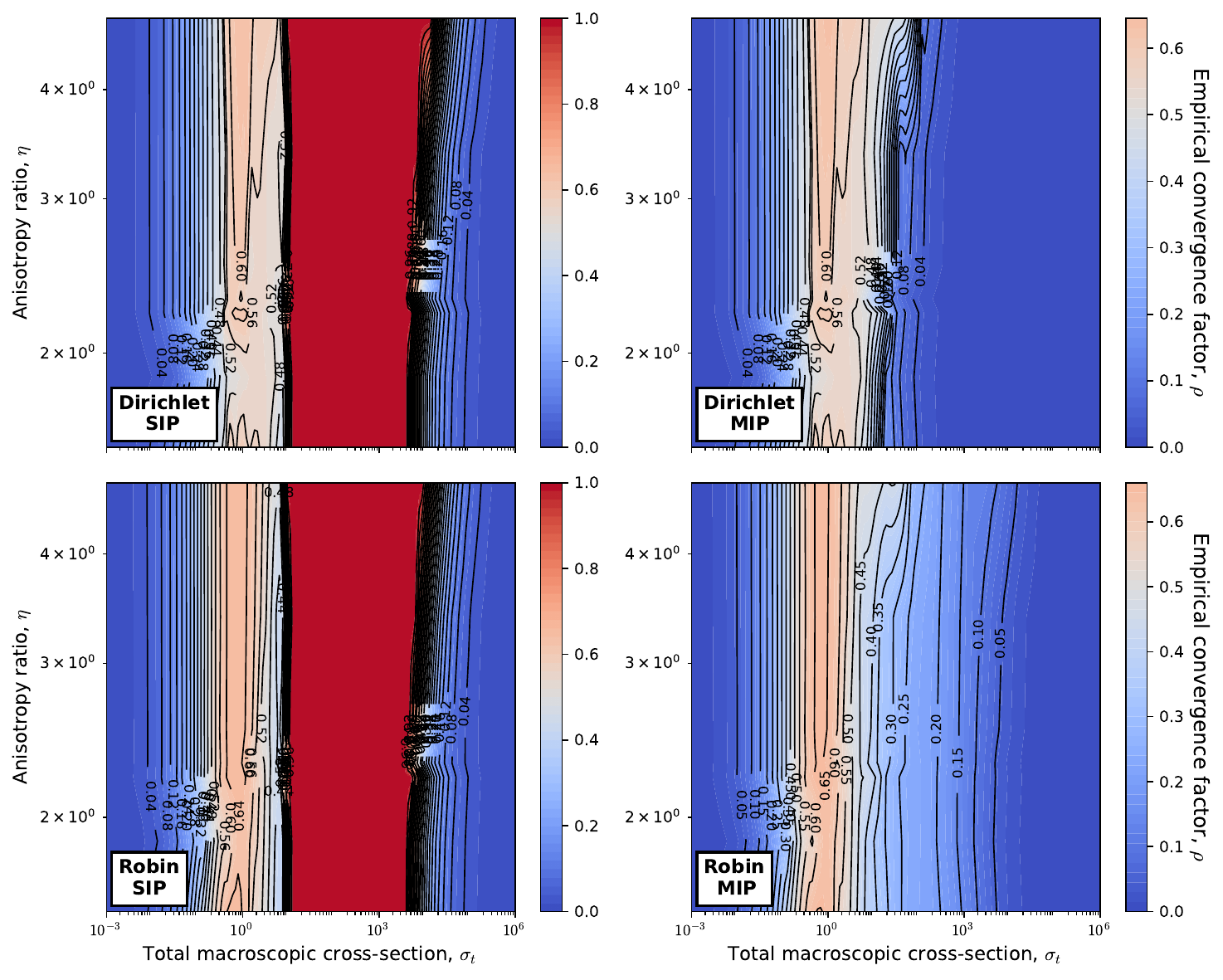}
  \caption{Experiment~\ref{subsec:exp7_anisotropy}. Empirical
    convergence factor as a function of $\sigma_t$, shown for multiple
    meshes with different anisotropy ratios $\eta$ (defined in
    \eqref{eqn:anisotropy_ratio}), comparing the four DSA variants
    (SIP/MIP combined with Dirichlet/Marshak diffusion boundary
    conditions).}
  \label{fig:exp7_a_comparison}
\end{figure}

To generate a controlled anisotropy sweep, we begin from a uniformly
random site distribution and generate bounded Voronoi tessellations
with a fixed element count $|\mathcal{T}_h|=64$. We then apply between
1 and 7 Lloyd iterations to obtain meshes spanning a broad range of
anisotropy values, while reducing the likelihood of numerical
artefacts arising from alignment between facets and discrete
ordinates. For all runs we fix $p=1$, $N_Q=16$ and $c=0.999$, and vary
over both $\sigma_t$ and $\eta$. The observed convergence factors for
the four DSA variants are shown in Figure~\ref{fig:exp7_a_comparison}.

Figure~\ref{fig:exp7_a_comparison} indicates that anisotropy has a
subtle but consequential influence on the convergence behaviour. As
$\eta$ increases, the stability region shrinks and convergence
degrades. The effect is most pronounced for the SIP variants, which
are divergent over a wider range of $\sigma_t$ as $\eta$
increases. The MIP variants mitigate this behaviour, although the
transition region widens with increasing $\eta$. In the
low-anisotropy regime, the behaviour is consistent with the results on
Lloyd-smoothed meshes reported earlier.

More specifically, in this experiment the SIP variants diverge around
$\sigma_t\sim 7$ and regain convergence once $\sigma_t/\eta \sim
2\times 10^3$. Prior to divergence, the maximal convergence factors
observed for the Dirichlet and Marshak SIP variants are $0.63$ and
$0.66$ respectively. For the MIP variants, the average convergence
factors are $0.14$ (Dirichlet) and $0.17$ (Marshak), with maximal
values of $0.63$ and $0.66$ in both cases. These maximal values remain broadly
consistent across the range of $\eta$, with a slight improvement in
the MIP convergence as $\eta$ increases.

Highly anisotropic meshes can also influence the global mesh diameter
$h$, so we additionally report trends using an $h$-independent
regularity indicator, the isoperimetric quotient. In the meshes used
elsewhere in this paper we typically observe isoperimetric ratios
above $0.65$. We find qualitatively similar behaviour under this
metric. As the isoperimetric quotient decreases (that is, as elements
become less regular), the range of $\sigma_t$ over which SIP-based
DSA diverges increases, while the MIP variants display comparable
robustness trends to those seen under the anisotropy ratio $\eta$.

Despite the known challenges associated with anisotropic
elements for both DG discretisations and DSA, these results indicate
that MIP-based DSA remains viable on significantly distorted Voronoi
meshes. More robust penalty constructions are likely to further
improve robustness in this setting; we refer to \cite{dsaanal} for
additional discussion. Finally, we note that diffusion solves can
become significantly more expensive on highly anisotropic meshes when
iterative solvers are used, due to increased conditioning. This is not
reflected here because we employ a direct solver, but it is an
important consideration for reproducibility and for large-scale
implementations.

\section{Conclusion and Discussion}
\label{section:conclusion}

We have studied diffusion synthetic acceleration (DSA) for
monoenergetic isotropic $S_N$ transport discretised by polytopal
discontinuous Galerkin methods. The focus has been on algorithmic
behaviour and implementation choices, with the accompanying analytic
results deferred to \cite{dsaanal}.

Across all experiments we compared four DSA variants, combining SIP or
MIP diffusion discretisations with either homogeneous Dirichlet or
Marshak (Robin) diffusion boundary conditions. The main conclusion is
that MIP-based DSA remains robust across the parameter ranges tested,
whereas SIP-based DSA exhibits a clear loss of robustness in the
intermediate regime and can diverge as optical thickness increases.

The numerical study covered variations in $\sigma_t$, the scattering
ratio $c$, the angular resolution $N_Q$, mesh refinement, polynomial
degree $p$ and mesh anisotropy. In regimes where source iteration
slows markedly, DSA reduces the observed convergence factor
substantially. The per-iteration diffusion overhead is largely
independent of $N_Q$, so its cost becomes increasingly amortised as
angular resolution increases.

Several broader trends also emerged. When plotted against optical
thickness $h\sigma_t$, the observed convergence factors display an
approximately invariant dependence across mesh refinements. MIP-based
DSA remains effective under $p$-refinement, with only mild variation
in convergence behaviour as $p$ increases. Strong mesh anisotropy
degrades performance and enlarges the instability region for SIP, but
MIP mitigates this effect and remains viable on distorted Voronoi
meshes. The choice of diffusion boundary condition has a secondary but
measurable effect. Using either Dirichlet or Marshak conditions can 
improve convergence in the intermediate regime depending on the context, 
while in the strongly diffusive limit the two boundary treatments become 
increasingly similar.

Future work will focus on analytic prediction of the observed
convergence behaviour, including sharper links between discrete upwind
transport fluxes and diffusion penalties, refined boundary layer
analysis for Marshak conditions on polytopal meshes and extension to
anisotropic scattering with diffusion coefficient
\[
D=\frac{1}{d\sigma_{\mathrm{tr}}},
\qquad
\sigma_{\mathrm{tr}}=\sigma_t-\sigma_{s,1}.
\]

\section{Data Availability}

A snapshot of the code used to generate all the data in this paper is
included in the Zenodo archive linked to this paper
\cite{zenodo}. This is part the Python package {\tt reyna}
\cite{reyna} for polytopic DGFEM.

\section{Acknowledgements}

ME is supported by scholarships from SAMBa under the project
EP/S022945/1. ME is also partially funded by the French Alternative
Energies and Atomic Energy Commission (CEA). TP received support from
the EPSRC programme grant EP/W026899/1, the Leverhulme RPG-2021-238
and EPSRC grant EP/X030067/1. All of this support is gratefully
acknowledged.

\bibliography{Literatur}
    
\end{document}